\documentclass[11pt]{amsart}

\usepackage{epigamath}

\usepackage[english]{babel}

\numberwithin{equation}{section}
\usepackage[mathscr]{eucal}
\usepackage{amssymb}

\usepackage{amsthm}% allows theorem* , e.g.
\usepackage{bbold}% allows \unit \mathbb{1}

\usepackage{enumitem}% allows resuming enumerate

\usepackage{amsmath}% allows pmatrix

\usepackage{tikz-cd}
\usetikzlibrary{decorations.pathreplacing,decorations.markings}
\usetikzlibrary{backgrounds}

\usepackage{tikzit}
\tikzstyle{bullet}=[fill=none, draw=black, shape=circle]
\tikzstyle{circle}=[fill=white, draw=black, shape=circle]
\tikzstyle{triangle}=[fill=white, draw=black, shape=regular polygon, regular polygon sides=3, shape border rotate=180, minimum width=4cm]
\tikzstyle{small triangle}=[fill=white, draw=black, shape=regular polygon, regular polygon sides=3, shape border rotate=180, minimum width=2cm]
\tikzstyle{tiny triangle}=[fill=white, draw=black, shape=regular polygon, regular polygon sides=3, shape border rotate=180, minimum width=1cm]
\tikzstyle{straight triangle}=[fill=white, draw=black, shape=regular polygon, regular polygon sides=3, minimum width=4cm]
\tikzstyle{straight small triangle}=[fill=white, draw=black, shape=regular polygon, regular polygon sides=3, minimum width=2cm]
\tikzstyle{straight tiny triangle}=[fill=white, draw=black, shape=regular polygon, regular polygon sides=3, minimum width=1cm]
\tikzstyle{tiny circle}=[fill=white, draw=black, shape=circle, scale=0.5]
\tikzstyle{really tiny circle}=[fill=white, draw=black, shape=circle, scale=0.25]
\tikzstyle{small rectangle}=[fill=white, draw=black, shape=rectangle, minimum height=1.5cm, minimum width=0.85cm]
\tikzstyle{tiny rectangle}=[fill=white, draw=black, shape=rectangle, minimum height=0.75cm, minimum width=0.70cm]

\tikzstyle{forward arrow}=[postaction={on each segment={mid farrow=black}}]
\tikzstyle{backward arrow}=[postaction={on each segment={mid barrow=black}}]
\tikzstyle{forward2 arrow}=[postaction={on each segment={mid farrow2=black}}]
\tikzstyle{backward2 arrow}=[postaction={on each segment={mid barrow2=black}}]
\tikzstyle{forward7 arrow}=[postaction={on each segment={mid farrow7=black}}]
\tikzstyle{backward7 arrow}=[postaction={on each segment={mid barrow7=black}}]

\tikzset{
	vertical align/.style={
	baseline=-.5*(height("$+$")-depth("$+$"))
	}
}
\tikzset{
    every picture/.style={vertical align}
}

\tikzset{
  on each segment/.style={
    decorate,
    decoration={
      show path construction,
      moveto code={},
      lineto code={
        \path [#1]
        (\tikzinputsegmentfirst) -- (\tikzinputsegmentlast);
      },
      curveto code={
        \path [#1] (\tikzinputsegmentfirst)
        .. controls
        (\tikzinputsegmentsupporta) and (\tikzinputsegmentsupportb)
        ..
        (\tikzinputsegmentlast);
      },
      closepath code={
        \path [#1]
        (\tikzinputsegmentfirst) -- (\tikzinputsegmentlast);
      },
    },
  },
  mid farrow2/.style={postaction={decorate,decoration={
        markings,
        mark=at position .25 with {\arrow[#1]{stealth}}
      }}},
  mid farrow7/.style={postaction={decorate,decoration={
        markings,
        mark=at position .75 with {\arrow[#1]{stealth}}
      }}},
  mid farrow/.style={postaction={decorate,decoration={
        markings,
        mark=at position .5 with {\arrow[#1]{stealth}}
      }}},
  mid barrow7/.style={postaction={decorate,decoration={
        markings,
        mark=at position .75 with {\arrowreversed[#1]{stealth}}
      }}},
  mid barrow2/.style={postaction={decorate,decoration={
        markings,
        mark=at position .25 with {\arrowreversed[#1]{stealth}}
      }}},
  mid barrow/.style={postaction={decorate,decoration={
        markings,
        mark=at position .5 with {\arrowreversed[#1]{stealth}}
      }}},
}

\tikzset{
    labelrotatebelow/.style={anchor=north, rotate=90, inner sep=1.0mm}
}
\tikzset{
    labelrotateabove/.style={anchor=south, rotate=90, inner sep=1.0mm}
}

\usepackage{graphicx}
\usepackage{mathtools}

\usepackage{extarrows}

\newtheorem{Thm}[equation]{Theorem}
\newtheorem*{Thm*}{Theorem}
\newtheorem{Prop}[equation]{Proposition}
\newtheorem{Lem}[equation]{Lemma}
\newtheorem{Cor}[equation]{Corollary}

\theoremstyle{definition}
\newtheorem{Def}[equation]{Definition}
\newtheorem{Ter}[equation]{Terminology}
\newtheorem{Not}[equation]{Notation}

\theoremstyle{remark}
\newtheorem{Rem}[equation]{Remark}

\newtheorem{Exa}[equation]{Example}

\newcommand{\nc}{\newcommand}
\nc{\dmo}{\DeclareMathOperator}

\renewcommand{\emptyset}{\varnothing}

\newcommand{\overbar}[1]{\mkern 1.5mu\overline{\mkern-1.5mu#1\mkern-1.5mu}\mkern 1.5mu}

\nc{\unitC}{\unit_{\cat C}}%
\nc{\unitD}{\unit_{\cat D}}%
\nc{\Pone}{{\mathbb{P}^1}}

\nc{\kapref}[1]{\hyperref[it:kappa#1]{\ensuremath{(\kappa#1)}}}
\nc{\sigref}[1]{\hyperref[it:sigma#1]{\ensuremath{(\sigma#1)}}}

\nc{\hooklongrightarrow}{\lhook\joinrel\longrightarrow}
\nc{\hooklongleftarrow}{\longleftarrow\joinrel\rhook}
\dmo{\gen}{gen}
\dmo{\StMod}{StMod}
\nc{\SpEn}{\cat S_{E(n)}}
\nc{\SpEnf}{\cat S_n}
\nc{\Lcomp}{L^{\mathrm{com}}} %I made this and the next one commands because I'm unsure of the choice of notation
\nc{\Ucomp}{U^{\mathrm{com}}}
\nc{\Loco}[1]{\Loc_{\otimes}\hspace{-0.3ex}\langle #1 \rangle}
\nc{\bbullet}{{\scriptscriptstyle\hspace{-1pt}\bullet}}
\nc{\bullett}{{\scriptscriptstyle\bullet}\hspace{-1pt}}
\nc{\LF}{L\hspace{-0.2ex}F}
\nc{\SpG}{\Sp^G}
\nc{\Prst}{{\cat P}\mathrm{r^{st}}}
\nc{\Mack}{\mathcal{M}ack}
\nc{\SC}{S\cat C}
\nc{\eps}{\epsilon}
\dmo{\lcm}{lcm}
\dmo{\tr}{tr}
\dmo{\Tr}{Tr}
\dmo{\fin}{{fin}}
\dmo{\DM}{DM}
\nc{\DMQ}{\DM_Q}
\dmo{\DerKal}{DMack}
\dmo{\Der}{D}
\dmo{\Derqc}{D_{qc}}
\dmo{\Derperf}{D_{perf}} %% not used 
\dmo{\DMot}{DMot}
\dmo{\rmH}{H}
\dmo{\piu}{\underline{\pi}}
\dmo{\Sphere}{\mathbb{S}}
\nc{\HA}{{\rmH \hspace{-0.2em}\bbA}}
\nc{\HZ}{{\rmH \hspace{-0.2em}\bbZ}}
\nc{\HZbar}{{\rmH \hspace{-0.2em}\underline{\bbZ}}}
\nc{\Fp}{{\bbF_{\hspace{-0.1em}p}}}
\nc{\HFp}{{\rmH \hspace{-0.15em}\bbF_{\hspace{-0.1em}p}}}
\nc{\DHZG}{\Der(\HZ_G)}
\nc{\DHZH}{\Der(\HZ_H)}
\nc{\DHZK}{\Der(\HZ_K)}
\nc{\DHZGN}{\Der(\HZ_{G/N})}
\nc{\DHZGG}{\Der(\HZ_{G/G})}
\nc{\DHZCp}{\Der(\HZ_{C_p})}
\nc{\DHZGprime}{\Der(\HZ_{G'})}
\nc{\DHZ}{\Der(\HZ)}
\nc{\frakp}{\mathfrak{p}}
\nc{\frakq}{\mathfrak{q}}
\nc{\Z}{\mathbb{Z}}
\nc{\SSG}{\text{sSet}_*^G}
\nc{\sSet}{\text{sSet}}

\dmo{\csupp}{csupp}
\dmo{\Con}{Conj}
\dmo{\Id}{Id}
\dmo{\Loc}{Loc}
\dmo{\rmK}{\textrm{\rm K}}
\dmo{\Spc}{Spc}
\dmo{\thick}{thick}
\nc{\thickt}[1]{\thick_\otimes\langle #1 \rangle}
\dmo{\cone}{cone}
\dmo{\End}{End}
\dmo{\Mor}{Mor}
\dmo{\Hom}{Hom}
\dmo{\id}{id}
\dmo{\incl}{incl}
\dmo{\Img}{Im}
\dmo{\im}{im}
\dmo{\Ker}{Ker}
\dmo{\ind}{ind}
\dmo{\CoInd}{coind}
\dmo{\res}{res}
\dmo{\infl}{infl}
\dmo{\triv}{triv}
\dmo{\Tel}{Tel} %telescope
\dmo{\grMod}{grMod}%
\dmo{\Mod}{Mod}%
\dmo{\Free}{Free}%
\dmo{\Kleisli}{Kleisli}%
\dmo{\opname}{op}
\dmo{\SH}{SH}% ground name for cat of spectra
\dmo{\smallb}{b}% ground exponent for ``bounded''
\dmo{\Spec}{Spec}
\dmo{\supp}{supp}
\dmo{\Supp}{Supp}
\nc{\SHA}{\SH{}^{\bbA^{1}}}% stab. motivic hom. cat of spec tra
\nc{\SHc}{{\SH^c}}
\nc{\SHp}{{\SH_{(p)}}}
\nc{\SHcp}{{\SH^c_{(p)}}}
\nc{\SHG}{\SH(G)}
\nc{\SHGp}{\SH(G)_{(p)}}
\nc{\SHGc}{\SHG^c}
\nc{\SHGcp}{\SHG^c_{(p)}}
\nc{\quadtext}[1]{\quad\textrm{#1}\quad}
\nc{\qquadtext}[1]{\qquad\textrm{#1}\qquad}
\nc{\adj}{\dashv}
\nc{\adjto}{\rightleftarrows}
\nc{\bbL}{\mathbb{L}}
\nc{\bbR}{\mathbb{R}}
\nc{\bbA}{\mathbb{A}}
\nc{\bbC}{\mathbb{C}}
\nc{\bbE}{\mathbb{E}}
\nc{\bbN}{\mathbb{N}}
\nc{\bbQ}{\mathbb{Q}}
\nc{\bbZ}{\mathbb{Z}}
\nc{\bbF}{\mathbb{F}}
\nc{\cat}[1]{\mathscr{#1}}%or: \nc{\cat}[1]{\mathcal{#1}}
\nc{\ie}{{\textit{i.e.}, }}
\nc{\into}{\mathop{\rightarrowtail}}
\nc{\inv}{^{-1}}
\nc{\isoto}{\mathop{\overset{\sim}\to}}
\nc{\isotoo}{\mathop{\overset{\sim}\too}}
\nc{\onto}{\mathop{\twoheadrightarrow}}
\nc{\too}{\mathop{\longrightarrow}\limits}
\nc{\mapstoo}{\longmapsto}
\nc{\adh}[1]{\overline{#1}}% adherence
\nc{\adhpt}[1]{\adh{\{#1\}}}% adherence of a pt
\nc{\aka}{{a.\,k.\,a.}\ }
\nc{\calF}{\mathcal{F}}
\nc{\eg}{{\sl e.\,g.}}
\nc{\Homcat}[1]{\Hom_{\cat #1}}
\nc{\hook}{\hookrightarrow}
\nc{\ideal}[1]{\langle #1\rangle}
\nc{\ihom}{{\underline{\hom}}}
\nc{\Mid}{\,\big|\,}
\nc{\KKleisli}{\,\text{-}\Kleisli}%
\nc{\FFree}{\,\text{-}\Free}%
\nc{\MMod}{\,\text{-}\Mod}%
\nc{\mmod}{\,\text{-mod}}%
\nc{\op}{^{\opname}}
\nc{\oto}[1]{\overset{#1}\to}
\nc{\otoo}[1]{\overset{#1}{\,\too\,}}
\nc{\sminus}{\!\smallsetminus\!}
\nc{\poplus}[1]{^{\oplus #1}}%
\nc{\potimes}[1]{^{\otimes #1}}% tensor power
\nc{\sbull}{{\scriptscriptstyle\bullet}}%\mathbf{\cdot}}%{}}
\nc{\SET}[2]{\big\{\,#1\Mid#2\,\big\}} 
\nc{\SpcK}{\Spc(\cat K)}% most used
\nc{\then}{\Rightarrow}
\nc{\unit}{\mathbb{1}}% unit for \otimes
\nc{\xra}{\xrightarrow}
\nc{\phigeom}[1]{\widetilde{\Phi}^{#1}}
\nc{\phigeomb}[1]{\Phi^{#1}}
\dmo{\Oname}{O}
\dmo{\proper}{proper}% for proper subgroups
\dmo{\lenormal}{\unlhd}
\dmo{\lnormal}{\lhd}
\nc{\normal}{\trianglelefteq}%\lhd
\nc{\Op}{\Oname^p}% O^p for maximal p-normal subgroup
\nc{\Oq}{\Oname^q}% as above for p=q
\dmo{\Sp}{Sp}
\dmo{\Ho}{Ho}
\dmo{\Fin}{Fin}
\dmo{\add}{add}
\dmo{\Fun}{Fun}
\dmo{\Ext}{Ext}
\dmo{\CAlg}{CAlg}
\dmo{\CMon}{CMon}
\dmo{\CC}{\cat C} %beren: I changed these, but left the O
\dmo{\DD}{\cat D}
\dmo{\OO}{\mathcal{O}}
\dmo{\Map}{Map}
\dmo{\Span}{Span}
\dmo{\N}{N}
\dmo{\Cat}{Cat}
\dmo{\colim}{colim}
\dmo{\hocolim}{hocolim}
\dmo{\Ch}{Ch}
\dmo{\A}{\mathbb{A}^{\mathrm{eff}}}
\nc{\AGeff}{\mathbb{A}_G^{\mathrm{eff}}}
\nc{\BGeff}{\mathcal{B}_G^{\mathrm{eff}}}
\nc{\BG}{{\mathcal{B}_G}}
\nc{\NBGeff}{{\N}{\BGeff}}
\dmo{\Ab}{Ab}
\dmo{\Set}{Set}
\dmo{\ev}{ev}
\dmo{\coev}{coev}
\dmo{\Spcl}{Spcl}
\nc{\Funadd}{\Fun_{\add}}
\dmo{\proj}{proj}
\dmo{\cof}{cof}
\def\isolow{\vbox to 0pt{\vss\hbox{$\scriptstyle\sim$}\vskip-2pt}}
\newcommand{\isor}{\xrightarrow{\isolow}}

\newcommand*{\xeq}[2][]{\underset{#1}{=}}

\newcommand{\supth}[1]{\ensuremath{#1^{\mathrm{th}}}}

\newcounter{enum-resume-hack}

\newlength{\mylength}
\newlength{\mylengthb}
\newlength{\mylengthc}
\newlength{\mylengthd}
\newlength{\mylengthe}
\newlength{\mylengthf}
\newlength{\mylengthg}
\newlength{\mylengthh}
\newlength{\mylengthi}
\newlength{\mylengthj}
\newlength{\mylengthk}
\EpigaVolumeYear{6}{2022} \EpigaArticleNr{18} \ReceivedOn{July 1,
2021}
\InFinalFormOn{May 31, 2022}
\AcceptedOn{June 21, 2022}

\title{A characterization of finite \'{e}tale morphisms\\ in tensor triangular geometry}
\titlemark{Finite \'{e}tale morphisms in tensor triangular geometry} 

\author{Beren Sanders}

\address{Mathematics Department, University of California, Santa Cruz, CA 95064, USA}
\email{beren@ucsc.edu}

\authormark{B.~Sanders}

\AbstractInEnglish{We provide a characterization of finite \'{e}tale morphisms in tensor triangular geometry.  They are precisely those functors which have a conservative right adjoint, satisfy Grothendieck--Neeman duality, and for which the relative dualizing object is trivial (via a canonically-defined map).}


\acknowledgement{The author is supported by NSF grant~DMS-1903429.}

\begin{document}

%%%%%%%%%%%%%%%%%%%%%%%%%%%%%%%
% Title page
%%%%%%%%%%%%%%%%%%%%%%%%%%%%%%%

\maketitle

\begin{prelims}

\DisplayAbstractInEnglish

\bigskip

\DisplayKeyWords

\medskip

\DisplayMSCclass

\end{prelims}

%%%%%%%%%%%%%%%%%%%%%
% Table of Contents
%%%%%%%%%%%%%%%%%%%%%

\newpage

\setcounter{tocdepth}{2}

\tableofcontents

%%%%%%%%%%%%%%%%%%%%%
% Content begins here
%%%%%%%%%%%%%%%%%%%%%

\section{Introduction}\label{sec:introduction} 

The purpose of this note is to give a characterization of ``finite \'{e}tale morphisms'' in tensor triangular geometry. We follow the notation, terminology, and perspective of \cite{BalmerDellAmbrogioSanders16}. In particular, we will work in the context of rigidly-compactly generated tensor-triangulated categories \cite[Definition~2.7]{BalmerDellAmbrogioSanders16}. The kind of characterization we have in mind is analogous to the following well-known characterization of smashing localizations:

\begin{Thm}\label{thm:smashing}
	Smashing localizations of a rigidly-compactly generated tensor-triangulated category $\cat T$ are precisely those geometric functors $f^*\colon\cat T \to \cat S$ between rigidly-compactly generated tensor-triangulated categories whose right adjoint $f_*$ is fully faithful.
\end{Thm}

We will recall a proof in Remark~\ref{rem:thm-1} below. Smashing localizations include, for example, restriction to a quasi-compact open subset of the Balmer spectrum. More generally, the tensor-triangular analogue of an \'{e}tale morphism is extension-of-scalars with respect to a commutative separable algebra (with smashing localizations being the special case of idempotent algebras). Finite \'{e}tale morphisms are, by definition, extension-of-scalars with respect to a \emph{compact} commutative separable algebra (see Definition~\ref{def:finite-etale}). Most smashing localizations are not \emph{finite} \'{e}tale morphisms, just as most open immersions are not proper. We will prove:

\begin{Thm}
	Finite \'{e}tale extensions of a rigidly-compactly generated tensor-triangulated category $\cat T$ are precisely those geometric functors $f^*\colon\cat T \to \cat S$ between rigidly-compactly generated tensor-triangulated categories which satisfy the following three properties:
	\begin{enumerate}
		\item \label{it:thm-a} $f^*$ satisfies Grothendieck--Neeman duality; 
		\item \label{it:thm-b} the right adjoint $f_*$ is conservative; 
		\item \label{it:thm-c} the canonical map $\unit_{\cat S}\to \omega_f$ is an isomorphism.
	\end{enumerate}
\end{Thm}

The terminology and notation will be explained in Section~\ref{sec:finite-etale}. We just remark that under hypothesis \hyperref[it:thm-a]{(a)}, the algebra $f_*(\unit_{\cat S})$ is rigid (a.k.a.~dualizable) and hence has an associated trace map. This corresponds by adjunction to a canonical map $\unit_{\cat S} \to \omega_f$ from the unit to the relative dualizing object, which hypothesis \hyperref[it:thm-c]{(c)} asserts is an isomorphism.

The keys to the theorem are the robust monadicity theorems which hold for triangulated categories and a deeper understanding of strongly separable algebras. Indeed, we begin the paper in Section~\ref{sec:strongly-separable} with a treatment of strongly separable algebras in arbitrary symmetric monoidal categories which may be of independent interest. We prove, in particular, that a rigid commutative algebra is separable if and only if it is strongly separable if and only if its canonically-defined trace form is nondegenerate (Corollary~\ref{cor:rigid-comm-alg}). We then turn in Section~\ref{sec:separable-triangulated} to tensor-triangulated categories and the role separable algebras play in that setting. A key tool is a strengthened version of separable monadicity (Proposition~\ref{prop:tensor-monadic}). We define finite \'{e}tale morphisms and prove the main theorem (Theorem~\ref{thm:main-thm}) in Section~\ref{sec:finite-etale}. We also show that if the target category is locally monogenic, then the conservativity condition \hyperref[it:thm-b]{(b)} is implied by the other two conditions (Corollary~\ref{cor:monogenic-char}). In Section~\ref{sec:examples}, we illustrate the theorem by giving some examples and non-examples of finite \'{e}tale morphisms in equivariant homotopy theory, algebraic geometry, and derived algebra.

\subsection*{Acknowledgments}
The author readily thanks Paul Balmer, Tobias Barthel, Ivo Dell'Ambrogio, Drew Heard, and Amnon Neeman. He is also grateful to Maxime Ramzi for pointing out an issue in the first proof of Proposition~\ref{prop:tensor-monadic} and to the referees for their helpful comments.

\section{Strongly separable algebras}\label{sec:strongly-separable}

We begin with a discussion of separable algebras in an arbitrary symmetric monoidal category. Although separable algebras are well-understood at this level of generality, we would like to clarify the notion of strongly separable algebra. Our main goal is to show that the equivalent characterizations of classical strongly separable algebras over fields established by \cite{Aguiar00} have suitable generalizations to arbitrary symmetric monoidal categories. The main punch line is that a rigid commutative algebra is separable if and only if it is strongly separable if and only if its trace form is nondegenerate (see Corollary~\ref{cor:rigid-comm-alg}). Moreover, this is the case if and only if it has the (necessarily unique) structure of a special symmetric Frobenius algebra.

\begin{Ter}
	Throughout this section, we work in a fixed symmetric monoidal category $(\cat C,\otimes,\unit)$. The symmetry isomorphism will be denoted 
	$\tau \colon A \otimes B \isor B \otimes A$. An object $A$ in $\cat C$ is \emph{rigid} (a.k.a.~\emph{dualizable}) if there exists an object $DA$ and morphisms $\eta\colon\unit \to DA \otimes A$ and $\epsilon\colon A \otimes DA \to \unit$ such that the composites
		\[
			A \xrightarrow{1 \otimes \eta} A\otimes DA \otimes A \xrightarrow{\epsilon\otimes 1} A
				\quad \text{ and }\quad
			DA \xrightarrow{\eta \otimes 1} DA\otimes A \otimes DA \xrightarrow{1 \otimes \epsilon} DA
		\]
	are the identity morphisms. It follows that the functor $DA \otimes -\colon\cat C \to \cat C$ is right adjoint to $A \otimes -\colon\cat C \to \cat C$. An \emph{algebra} $A$ is an associative unital monoid in~$\cat C$. The multiplication and unit maps will be denoted 
	$\mu\colon A \otimes A \to A$ and $u\colon \unit \to A$.
\end{Ter}

\begin{Def}
	An algebra $(A,\mu,u)$ is \emph{separable} if there exists a map $\sigma\colon A\to A \otimes A$ such that 
		\begin{enumerate}
			\item[($\sigma$1)]\label{it:sigma1} $\mu \circ \sigma = \id_A$ and
			\item[($\sigma$2)]\label{it:sigma2} $(1\otimes \mu)\circ(\sigma\otimes 1) = \sigma \circ \mu = (\mu \otimes 1)\circ(1\otimes \sigma)$ as maps $A\otimes A \to A \otimes A$.
		\end{enumerate}
	In other words, $A$ is separable if the multiplication map $\mu\colon A\otimes A\to A$ admits an $(A,A)$-bilinear section.
\end{Def}

\begin{Rem}\label{rem:precompose}
	If we precompose such a section $\sigma$ with the unit $u\colon \unit \to A$, we obtain a map \[\kappa\coloneqq \sigma\circ u \colon  \unit \to A \otimes A\] which satisfies
		\begin{enumerate}
			\item[($\kappa$1)]\label{it:kappa1} $\mu \circ\kappa = u$, and
			\item[($\kappa$2)]\label{it:kappa2} $(1\otimes \mu)\circ(\kappa\otimes 1) = (\mu \otimes 1)\circ(1\otimes \kappa)$ as maps $A \to A \otimes A$.
		\end{enumerate}
	Conversely, given such a $\kappa$, the map $A \to A \otimes A$ displayed in \kapref2 satisfies axioms \sigref1 and \sigref2.  Thus, an algebra $A$ is separable if and only if it admits a map ${\kappa \colon  \unit \to A \otimes A}$ satisfying \kapref1 and \kapref2.  Such a map $\kappa$ is called a \emph{separability idempotent}.
\end{Rem}

\begin{Rem}
	We refer the reader to \cite{AuslanderGoldman60, ChaseHarrisonRosenberg65, KnusOjanguren74}, \cite[Chapter~10]{Pierce82}, and \cite{Ford17} for further information about separable algebras and their role in classical representation theory and algebraic geometry.  The following notion of a \emph{strongly} separable algebra was originally studied by Kanzaki and Hattori \cite{Hattori65,Kanzaki62}: 
\end{Rem}

\begin{Def}
	An algebra $A$ is \emph{strongly separable} if there exists a map $\kappa\colon \unit \to A \otimes A$ satisfying \kapref1, \kapref2,  and 
		\begin{enumerate}
			\item[($\kappa$3)]\label{it:kappa3} $\kappa = \tau \circ \kappa$.
		\end{enumerate}
	In other words, $A$ is strongly separable if it admits a \emph{symmetric} separability idempotent.
\end{Def}

\begin{Rem}
	For classical algebras over a field, \cite{Aguiar00} provides several equivalent characterizations of strongly separable algebras.  Our present goal is to clarify the extent to which these characterizations hold in an arbitrary symmetric monoidal category.  For this purpose, the graphical calculus of string diagrams will be very convenient.  We refer the reader to \cite[Sections~3--4]{Selinger11} and \cite[Section~2]{PontoShulman13} for more information concerning these diagrams and suffice ourselves to remark that it is a routine exercise to convert a proof involving string diagrams into a detailed proof using commutative diagrams.
\end{Rem}

\begin{Not}\label{not:string-diagrams}
	We will read our string diagrams from bottom to top.  The multiplication map $\mu\colon A \otimes A \to A$ and the map $\kappa\colon \unit \to A \otimes A$ will be represented by
		\[
		%\resizebox{!}{1.0cm}{%
		\begin{tikzpicture}[scale=0.5,baseline=(current bounding box.center)]
			\begin{pgfonlayer}{nodelayer}
				\node [style=none] (1) at (1, 3) {};
				\node [style=bullet] (4) at (2.5, 4.5) {$\mu$};
				\node [style=none] (6) at (4, 3) {};
				\node [style=none] (13) at (2.5, 6) {};
			\end{pgfonlayer}
			\begin{pgfonlayer}{edgelayer}
				\draw [in=180, out=90] (1.center) to (4);
				\draw [in=90, out=0] (4) to (6);
				%\draw (5.center) to (0.center);
				%\draw [in=90, out=-180, looseness=1.25] (6) to (7.center);
				%\draw [in=90, out=0, looseness=1.25] (6) to (8.center);
				%\draw [in=90, out=-90] (7.center) to (10.center);
				%\draw [in=90, out=-90] (8.center) to (9.center);
				\draw (13.center) to (4);
			\end{pgfonlayer}
		\end{tikzpicture}%}
		\quad\text{ and }\quad
		\begin{tikzpicture}[scale=0.25,baseline=(current bounding box.center)]
			\begin{pgfonlayer}{nodelayer}
				\node [fill=white, draw=black, shape=regular polygon, regular polygon sides=3, shape border rotate=180] (0) at (0, 0) {$\kappa$};
		%		\node [style=really tiny circle] (1) at (0, -2.5) {};
				\node [style=none] (2) at (-1, 4) {};
				\node [style=none] (3) at (1, 4) {};
				\node [style=none] (4) at (-1, 1) {};
				\node [style=none] (5) at (1, 1) {};
			\end{pgfonlayer}
			\begin{pgfonlayer}{edgelayer}
				\draw (4.center) to (2.center);
				\draw (5.center) to (3.center);
			\end{pgfonlayer}
		\end{tikzpicture}%}
		\]
	while the unit $u \colon \unit \to A$  and the identity $\id\colon A \to A$ will be represented by
		\[
			u = 
			\begin{tikzpicture}[scale=0.20,baseline=(current bounding box.center)]
				\begin{pgfonlayer}{nodelayer}
					\node [style=really tiny circle] (1) at (0, -2.5) {};
					\node [style=none] (2) at (0, 4) {};
				\end{pgfonlayer}
				\begin{pgfonlayer}{edgelayer}
					\draw (1.center) to (2.center);
				\end{pgfonlayer}
			\end{tikzpicture}
			\quad\text{ and }\quad
			\id_A =
			\begin{tikzpicture}[scale=0.20,baseline=(current bounding box.center)]
				\begin{pgfonlayer}{nodelayer}
					\node [style=none] (1) at (0, -2.5) {};
					\node [style=none] (2) at (0, 4) {};
				\end{pgfonlayer}
				\begin{pgfonlayer}{edgelayer}
					\draw (1.center) to (2.center);
				\end{pgfonlayer}
			\end{tikzpicture}
		\]
	Thus, for example, axiom \kapref2 reads
		\begin{equation}\label{eq:string-k2}
			\resizebox{!}{0.8cm}{
			\begin{tikzpicture}[baseline=(current bounding box.center)]
				\begin{pgfonlayer}{nodelayer}
		\node [style=none] (0) at (-1, 1) {};
		\node [style=none] (1) at (1, 1) {};
		\node [style=triangle] (2) at (0, 0) {};
%		\node [style=tiny circle] (3) at (0, -2) {};
		\node [style=bullet] (4) at (2.5, 2.5) {};
		\node [style=none] (5) at (-1, 4) {};
		\node [style=none] (6) at (4, -2) {};
		\node [style=none] (7) at (4, 1) {};
		\node [style=none] (8) at (2.5, 4) {};
	\end{pgfonlayer}
	\begin{pgfonlayer}{edgelayer}
		\draw [in=180, out=90] (1.center) to (4);
		\draw [in=90, out=0] (4) to (7.center);
		\draw (7.center) to (6);
		\draw (5.center) to (0.center);
		\draw (8.center) to (4);
	\end{pgfonlayer}
%
%\begin{pgfonlayer}{nodelayer}
%		\node [style=none] (0) at (-1, 1) {};
%		\node [style=none] (1) at (1, 1) {};
%		\node [style=triangle] (2) at (0, 0) {};
%		\node [style=tiny circle] (3) at (0, -2) {};
%		\node [style=bullet] (4) at (2.5, 2.5) {};
%		\node [style=none] (5) at (-1, 4) {};
%		\node [style=tiny circle] (6) at (4, -2) {};
%		\node [style=none] (7) at (4, 1) {};
%		\node [style=none] (8) at (2.5, 4) {};
%	\end{pgfonlayer}
%	\begin{pgfonlayer}{edgelayer}
%		\draw [in=180, out=90] (1.center) to (4);
%		\draw [in=90, out=0] (4) to (7.center);
%		\draw (7.center) to (6);
%		\draw (5.center) to (0.center);
%		\draw (8.center) to (4);
%	\end{pgfonlayer}

			\end{tikzpicture}}
			=
			\resizebox{!}{0.8cm}{
			\begin{tikzpicture}[baseline=(current bounding box.center)]
				\begin{pgfonlayer}{nodelayer}
		\node [style=none] (0) at (1, 1) {};
		\node [style=none] (1) at (-1, 1) {};
		\node [style=triangle] (2) at (0, 0) {};
%		\node [style=tiny circle] (3) at (0, -2) {};
		\node [style=bullet] (4) at (-2.5, 2.5) {};
		\node [style=none] (5) at (1, 4) {};
		\node [style=none] (6) at (-4, -2) {};
		\node [style=none] (7) at (-4, 1) {};
		\node [style=none] (8) at (-2.5, 4) {};
	\end{pgfonlayer}
	\begin{pgfonlayer}{edgelayer}
		\draw [in=0, out=90] (1.center) to (4);
		\draw [in=90, out=180] (4) to (7.center);
		\draw (7.center) to (6);
		\draw (5.center) to (0.center);
		\draw (8.center) to (4);
	\end{pgfonlayer}

			\end{tikzpicture}}
		\end{equation}
\end{Not}

\begin{Prop}\label{prop:strongly-sep-alt-axioms}
	An algebra $(A,\mu,u)$ is strongly separable if and only if there exists a morphism $\kappa\colon \unit \to A \otimes A$ satisfying \kapref2 and
		\begin{enumerate}
			\item[($\kappa$\textup{4})]\label{it:kappa4} $\mu \circ \tau \circ \kappa = u$.
		\end{enumerate}
\end{Prop}

\begin{proof}
	By definition, an algebra is strongly separable if it admits a morphism $\kappa$ satisfying \kapref1, \kapref2, and \kapref3.  It is immediate that \kapref1 and \kapref3 together imply~\kapref4.  It is also immediate that \kapref3 and \kapref4 together imply \kapref1.  Thus, the claim will be established if we can prove that \kapref2 and \kapref4 together imply~\kapref3.  Using Notation~\ref{not:string-diagrams}, observe:
		\begin{equation}\label{eq:string-row1}
			\begin{tikzpicture}[scale=0.20,transform shape]
				\begin{pgfonlayer}{nodelayer}
		\node [style=triangle] (0) at (0, 0) {};
	%	\node [style=tiny circle] (1) at (0, -2) {};
		\node [style=none] (2) at (-1, 4) {};
		\node [style=none] (3) at (1, 4) {};
		\node [style=none] (4) at (-1, 1) {};
		\node [style=none] (5) at (1, 1) {};
	\end{pgfonlayer}
	\begin{pgfonlayer}{edgelayer}
		\draw (4.center) to (2.center);
		\draw (5.center) to (3.center);
	\end{pgfonlayer}

			\end{tikzpicture}
			\xeq[(\text{unital})]{}
			\begin{tikzpicture}[scale=0.20,transform shape]
			\begin{pgfonlayer}{nodelayer}
		\node [style=none] (0) at (-1, 1) {};
		\node [style=none] (1) at (1, 1) {};
		\node [style=triangle] (2) at (0, 0) {};
%		\node [style=tiny circle] (3) at (0, -2) {};
		\node [style=bullet] (4) at (2.5, 2.5) {};
		\node [style=none] (5) at (-1, 4) {};
		\node [style=tiny circle] (6) at (4, -2) {};
		\node [style=none] (7) at (4, 1) {};
		\node [style=none] (8) at (2.5, 4) {};
	\end{pgfonlayer}
	\begin{pgfonlayer}{edgelayer}
		\draw [in=180, out=90] (1.center) to (4);
		\draw [in=90, out=0] (4) to (7.center);
		\draw (7.center) to (6);
		\draw (5.center) to (0.center);
		\draw (8.center) to (4);
	\end{pgfonlayer}
%
%\begin{pgfonlayer}{nodelayer}
%		\node [style=none] (0) at (-1, 1) {};
%		\node [style=none] (1) at (1, 1) {};
%		\node [style=triangle] (2) at (0, 0) {};
%		\node [style=tiny circle] (3) at (0, -2) {};
%		\node [style=bullet] (4) at (2.5, 2.5) {};
%		\node [style=none] (5) at (-1, 4) {};
%		\node [style=tiny circle] (6) at (4, -2) {};
%		\node [style=none] (7) at (4, 1) {};
%		\node [style=none] (8) at (2.5, 4) {};
%	\end{pgfonlayer}
%	\begin{pgfonlayer}{edgelayer}
%		\draw [in=180, out=90] (1.center) to (4);
%		\draw [in=90, out=0] (4) to (7.center);
%		\draw (7.center) to (6);
%		\draw (5.center) to (0.center);
%		\draw (8.center) to (4);
%	\end{pgfonlayer}

			\end{tikzpicture}
			\xeq[\kapref4]{}
			\begin{tikzpicture}[scale=0.20,transform shape]
				\begin{pgfonlayer}{nodelayer}
		\node [style=none] (0) at (-1, 3) {};
		\node [style=none] (1) at (1, 3) {};
		\node [style=triangle] (2) at (0, 2) {};
	%	\node [style=tiny circle] (3) at (0, 0) {};
		\node [style=bullet] (4) at (2.5, 4.5) {};
		\node [style=none] (5) at (-1, 6) {};
		\node [style=bullet] (6) at (4, 3) {};
		\node [style=none] (7) at (2.75, 1.5) {};
		\node [style=none] (8) at (5.25, 1.5) {};
		\node [style=none] (9) at (2.75, -1) {};
		\node [style=none] (10) at (5.25, -1) {};
		\node [style=triangle] (11) at (4, -2) {};
	%	\node [style=tiny circle] (12) at (4, -4) {};
		\node [style=none] (13) at (2.5, 6) {};
	\end{pgfonlayer}
	\begin{pgfonlayer}{edgelayer}
		\draw [in=180, out=90] (1.center) to (4);
		\draw [in=90, out=0] (4) to (6);
		\draw (5.center) to (0.center);
		\draw [in=90, out=-180, looseness=1.25] (6) to (7.center);
		\draw [in=90, out=0, looseness=1.25] (6) to (8.center);
		\draw [in=90, out=-90] (7.center) to (10.center);
		\draw [in=90, out=-90] (8.center) to (9.center);
		\draw (13.center) to (4);
	\end{pgfonlayer}

			\end{tikzpicture}
			\xeq[(\text{assoc})]{}
			\begin{tikzpicture}[scale=0.20,transform shape]
				\begin{pgfonlayer}{nodelayer}
		\node [style=none] (0) at (-3.25, -0.5) {};
		\node [style=none] (1) at (-1.25, -0.5) {};
		\node [style=triangle] (2) at (-2.25, -1.5) {};
	%	\node [style=tiny circle] (3) at (-2.25, -3.5) {};
		\node [style=bullet] (4) at (0.25, 1) {};
		\node [style=none] (5) at (-3.25, 4) {};
		\node [style=bullet] (6) at (1.75, 2.5) {};
		\node [style=none] (7) at (1.25, 0) {};
		\node [style=none] (8) at (3.5, 0) {};
		\node [style=none] (9) at (1.25, -3) {};
		\node [style=none] (10) at (3.75, -3) {};
		\node [style=triangle] (11) at (2.5, -4) {};
	%	\node [style=tiny circle] (12) at (2.5, -6) {};
		\node [style=none] (13) at (1.75, 4) {};
	\end{pgfonlayer}
	\begin{pgfonlayer}{edgelayer}
		\draw [in=180, out=90] (1.center) to (4);
		\draw (5.center) to (0.center);
		\draw [in=90, out=0, looseness=1.25] (6) to (8.center);
		\draw [in=90, out=-90] (7.center) to (10.center);
		\draw [in=90, out=-90] (8.center) to (9.center);
		\draw [in=90, out=0] (4) to (7.center);
		\draw [in=-180, out=90, looseness=1.25] (4) to (6);
		\draw (13.center) to (6);
	\end{pgfonlayer}

			\end{tikzpicture}
			\xeq[\kapref2]{}
			\begin{tikzpicture}[scale=0.20,transform shape]
				\begin{pgfonlayer}{nodelayer}
		\node [style=none] (0) at (-3.25, -0.5) {};
		\node [style=none] (1) at (-1.25, -0.5) {};
		\node [style=triangle] (2) at (-2.25, -1.5) {};
	%	\node [style=tiny circle] (3) at (-2.25, -3.5) {};
		\node [style=bullet] (4) at (0.25, 3) {};
		\node [style=bullet] (5) at (-4.25, 3) {};
		\node [style=none] (6) at (1.5, 0) {};
		\node [style=none] (7) at (-5.75, 1.75) {};
		\node [style=none] (8) at (3.75, 0) {};
		\node [style=none] (9) at (1.5, -3) {};
		\node [style=triangle] (10) at (2.75, -4) {};
	%	\node [style=tiny circle] (11) at (2.75, -6) {};
		\node [style=none] (12) at (0.25, 4) {};
		\node [style=none] (13) at (-4.25, 4) {};
		\node [style=none] (14) at (3.75, -3) {};
	\end{pgfonlayer}
	\begin{pgfonlayer}{edgelayer}
		\draw [in=180, out=90] (1.center) to (4);
		\draw [in=90, out=0, looseness=0.75] (5) to (0.center);
		\draw [in=90, out=-90] (6.center) to (9.center);
		\draw [in=90, out=-90, looseness=0.50] (7.center) to (8.center);
		\draw [in=90, out=0] (4) to (6.center);
		\draw (13.center) to (5);
		\draw (12.center) to (4);
		\draw [in=90, out=-180, looseness=1.25] (5) to (7.center);
		\draw (8.center) to (14.center);
	\end{pgfonlayer}

			\end{tikzpicture}
		\end{equation}
	We can then rearrange the last diagram by pulling the left-hand multiplication to the right-hand side and continue:
		\begin{equation}\label{eq:string-row2}
			\begin{tikzpicture}[scale=0.20,transform shape]
				\begin{pgfonlayer}{nodelayer}
		\node [style=none] (0) at (-4.25, -3) {};
		\node [style=none] (1) at (-2.25, -3) {};
		\node [style=triangle] (2) at (-3.25, -4) {};
%		\node [style=tiny circle] (3) at (-3.25, -6) {};
		\node [style=bullet] (4) at (-0.5, 0.5) {};
		\node [style=bullet] (5) at (5, 0.5) {};
		\node [style=none] (6) at (3.5, -3) {};
		\node [style=none] (7) at (1, -3) {};
		\node [style=triangle] (8) at (2.25, -4) {};
%		\node [style=tiny circle] (9) at (2.25, -6) {};
		\node [style=none] (10) at (3.25, 3) {};
		\node [style=none] (11) at (1.5, 3) {};
		\node [style=none] (12) at (6.25, -1) {};
		\node [style=none] (13) at (1.5, 4) {};
		\node [style=none] (14) at (3.25, 4) {};
	\end{pgfonlayer}
	\begin{pgfonlayer}{edgelayer}
		\draw [in=180, out=90] (1.center) to (4);
		\draw [in=90, out=-90, looseness=0.50] (11.center) to (5);
		\draw [in=90, out=-90, looseness=0.50] (10.center) to (4);
		\draw [in=90, out=-180, looseness=1.25] (5) to (6.center);
		\draw [in=90, out=0] (4) to (7.center);
		\draw [in=90, out=0, looseness=1.25] (5) to (12.center);
		\draw [in=-90, out=90, looseness=0.50] (0.center) to (12.center);
		\draw (13.center) to (11.center);
		\draw (14.center) to (10.center);
	\end{pgfonlayer}

			\end{tikzpicture}
			\xeq[\kapref2]{}
			\begin{tikzpicture}[scale=0.20,transform shape]
			\begin{pgfonlayer}{nodelayer}
		\node [style=none] (0) at (3.75, -0.75) {};
		\node [style=none] (1) at (1.75, -0.75) {};
		\node [style=triangle] (2) at (2.75, -1.75) {};
%		\node [style=tiny circle] (3) at (2.75, -3.75) {};
		\node [style=bullet] (4) at (-1.75, 2.5) {};
		\node [style=none] (5) at (3.75, 3.5) {};
		\node [style=triangle] (6) at (-3, -4.25) {};
%		\node [style=tiny circle] (7) at (-3, -6.25) {};
		\node [style=none] (8) at (-1.75, 3.5) {};
		\node [style=none] (9) at (3, 5.5) {};
		\node [style=none] (10) at (-0.5, 5.5) {};
		\node [style=none] (11) at (3, 6.25) {};
		\node [style=none] (12) at (-0.5, 6.25) {};
		\node [style=bullet] (13) at (0, 1) {};
		\node [style=none] (14) at (-4.25, -0.75) {};
		\node [style=none] (15) at (-1.75, -0.75) {};
		\node [style=none] (16) at (-4.25, -3.25) {};
		\node [style=none] (17) at (-1.75, -3.25) {};
		\node [style=none] (18) at (-4.25, 1) {};
	\end{pgfonlayer}
	\begin{pgfonlayer}{edgelayer}
		\draw (5.center) to (0.center);
		\draw (8.center) to (4);
		\draw [in=90, out=-90, looseness=0.50] (10.center) to (5.center);
		\draw [in=90, out=-90, looseness=0.50] (9.center) to (8.center);
		\draw (12.center) to (10.center);
		\draw (11.center) to (9.center);
		\draw [in=90, out=0, looseness=1.25] (4) to (13);
		\draw [in=90, out=0, looseness=1.25] (13) to (1.center);
		\draw [in=90, out=-90] (14.center) to (17.center);
		\draw [in=90, out=-90] (15.center) to (16.center);
		\draw [in=-180, out=90, looseness=1.25] (15.center) to (13);
		\draw (14.center) to (18.center);
		\draw [in=180, out=90] (18.center) to (4);
	\end{pgfonlayer}

			\end{tikzpicture}
			\xeq[(\text{assoc})]{}
			\begin{tikzpicture}[scale=0.20,transform shape]
			\begin{pgfonlayer}{nodelayer}
		\node [style=none] (0) at (3.5, 1) {};
		\node [style=none] (1) at (1.5, 1) {};
		\node [style=triangle] (2) at (2.5, 0) {};
%		\node [style=tiny circle] (3) at (2.5, -2) {};
		\node [style=bullet] (4) at (0, 2.5) {};
		\node [style=none] (5) at (3.5, 3.5) {};
		\node [style=bullet] (6) at (-1.5, 1) {};
		\node [style=none] (7) at (-2.75, -0.5) {};
		\node [style=none] (8) at (-0.25, -0.5) {};
		\node [style=none] (9) at (-2.75, -3) {};
		\node [style=none] (10) at (-0.25, -3) {};
		\node [style=triangle] (11) at (-1.5, -4) {};
%		\node [style=tiny circle] (12) at (-1.5, -6) {};
		\node [style=none] (13) at (0, 3.5) {};
		\node [style=none] (14) at (3.5, 5.5) {};
		\node [style=none] (15) at (0, 5.5) {};
		\node [style=none] (16) at (3.5, 6.25) {};
		\node [style=none] (17) at (0, 6.25) {};
	\end{pgfonlayer}
	\begin{pgfonlayer}{edgelayer}
		\draw [in=0, out=90] (1.center) to (4);
		\draw [in=90, out=-180] (4) to (6);
		\draw (5.center) to (0.center);
		\draw [in=90, out=-180, looseness=1.25] (6) to (7.center);
		\draw [in=90, out=0, looseness=1.25] (6) to (8.center);
		\draw [in=90, out=-90] (7.center) to (10.center);
		\draw [in=90, out=-90] (8.center) to (9.center);
		\draw (13.center) to (4);
		\draw [in=90, out=-90, looseness=0.50] (15.center) to (5.center);
		\draw [in=90, out=-90, looseness=0.50] (14.center) to (13.center);
		\draw (17.center) to (15.center);
		\draw (16.center) to (14.center);
	\end{pgfonlayer}

			\end{tikzpicture}
			\xeq[\kapref4]{}
			\begin{tikzpicture}[scale=0.20,transform shape]
				\begin{pgfonlayer}{nodelayer}
		\node [style=none] (0) at (2, 1) {};
		\node [style=none] (1) at (0, 1) {};
		\node [style=triangle] (2) at (1, 0) {};
		\node [style=bullet] (4) at (-1.5, 2.5) {};
		\node [style=none] (5) at (2, 3.5) {};
		\node [style=none] (6) at (-3, 1) {};
		\node [style=none] (8) at (-1.5, 3.5) {};
		\node [style=none] (9) at (2, 5.5) {};
		\node [style=none] (10) at (-1.5, 5.5) {};
		\node [style=none] (11) at (2, 6.25) {};
		\node [style=none] (12) at (-1.5, 6.25) {};
		\node [style=tiny circle] (13) at (-3, -1) {};
	\end{pgfonlayer}
	\begin{pgfonlayer}{edgelayer}
		\draw [in=0, out=90] (1.center) to (4);
		\draw [in=90, out=-180] (4) to (6.center);
		\draw (5.center) to (0.center);
		\draw (8.center) to (4);
		\draw [in=90, out=-90, looseness=0.50] (10.center) to (5.center);
		\draw [in=90, out=-90, looseness=0.50] (9.center) to (8.center);
		\draw (12.center) to (10.center);
		\draw (11.center) to (9.center);
		\draw (6.center) to (13.center);
	\end{pgfonlayer}

			\end{tikzpicture}
			\xeq[(\text{unital})]{}
			\begin{tikzpicture}[scale=0.20,transform shape]
			\begin{pgfonlayer}{nodelayer}
		\node [style=none] (0) at (1, 1) {};
		\node [style=none] (1) at (-1, 1) {};
		\node [style=triangle] (2) at (0, 0) {};
%		\node [style=tiny circle] (3) at (0, -2) {};
		\node [style=none] (4) at (1, 3.25) {};
		\node [style=none] (5) at (-1, 3.25) {};
		\node [style=none] (6) at (1, 5.5) {};
		\node [style=none] (7) at (-1, 5.5) {};
		\node [style=none] (8) at (1, 6.25) {};
		\node [style=none] (9) at (-1, 6.25) {};
	\end{pgfonlayer}
	\begin{pgfonlayer}{edgelayer}
		\draw (4.center) to (0.center);
		\draw [in=90, out=-90, looseness=0.75] (7.center) to (4.center);
		\draw [in=90, out=-90, looseness=0.75] (6.center) to (5.center);
		\draw (9.center) to (7.center);
		\draw (8.center) to (6.center);
		\draw (5.center) to (1.center);
	\end{pgfonlayer}

			\end{tikzpicture}
		\end{equation}
	This establishes $\kappa = \tau\circ \kappa$, which is axiom \kapref3.
\end{proof}

\begin{Cor}\label{cor:commutative-separable-is-strongly-separable}
	 Any commutative separable algebra is strongly separable.
\end{Cor}

\begin{proof}
	This follows from Proposition~\ref{prop:strongly-sep-alt-axioms} since axiom \kapref4 coincides with axiom~\kapref1 when the algebra is commutative.
\end{proof}

\begin{Rem}
	For string diagrams involving a rigid object $A$, we will use the direction of a string to indicate whether it represents $A$ or its dual $DA$.  For example, the unit $\unit \to DA \otimes A$ and counit $A \otimes DA \to \unit$ are represented by 
		\[
			\input{fig-dual-unitb.tikz}
			\quad\text{ and }\quad
			\input{fig-dual-counitb.tikz}
		\]
	respectively, and the unit-counit relations are given by
		
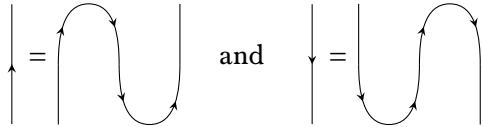
\[
			%\resizebox{!}{1.0cm}{%
			%\begin{tikzpicture}[framed, background rectangle/.style={draw=white},scale=0.5,baseline=(current bounding box.center)]
			\begin{tikzpicture}[scale=0.4,baseline=(current bounding box.center)]
				\begin{pgfonlayer}{nodelayer}
		\node  (0) at (0, 0) {};
		\node  (9) at (0, 4.65) {};
	\end{pgfonlayer}
	\begin{pgfonlayer}{edgelayer}
		\draw [style=forward arrow] (0) to (9);
	\end{pgfonlayer}

			\end{tikzpicture}
			%}
			=
			\,
			%\resizebox{!}{1.0cm}{%
			%\begin{tikzpicture}[framed, background rectangle/.style={double,ultra thick,draw=red, top color=blue, rounded corners},scale=0.5,baseline=(current bounding box.center)]
			%\begin{tikzpicture}[framed, background rectangle/.style={draw=white},scale=0.5,baseline=(current bounding box.center)]
			\begin{tikzpicture}[scale=0.4,baseline=(current bounding box.center)]
				\begin{pgfonlayer}{nodelayer}
		\node [style=none] (0) at (0, 0) {};
		\node [style=none] (1) at (-1, 2) {};
		\node [style=none] (2) at (1, 2) {};
		\node [style=none] (4) at (1, 2) {};
		\node [style=none] (8) at (1, 4) {};
		\node [style=none] (9) at (-2, 4) {};
		\node [style=none] (10) at (-3, 2) {};
		\node [style=none] (11) at (-1, 2) {};
		\node [style=none] (12) at (-1, 2) {};
		\node [style=none] (13) at (-3, 0) {};
	\end{pgfonlayer}
	\begin{pgfonlayer}{edgelayer}
		\draw [style=forward arrow] (1.center)
				[in=180, out=-90] to (0.center)
				[in=-90, out=0] to (2.center);
%		\draw [style=forward arrow, in=180, out=-90] (1.center) to (0);
%		\draw [style=forward arrow, in=-90, out=0] (0) to (2.center);
		\draw (8.center) to (4.center);
		\draw [style=forward arrow] (10.center)
				[in=-180, out=90] to (9.center)
				[in=90, out=0] to (11.center);
%		\draw [style=forward arrow, in=-180, out=90] (10.center) to (9);
%		\draw [style=forward arrow, in=90, out=0] (9) to (11.center);
		\draw (10.center) to (13.center);
	\end{pgfonlayer}

			\end{tikzpicture}
			%}
			\quad\text{ and }\quad
			%\resizebox{!}{1.0cm}{%
			\begin{tikzpicture}[scale=0.4,baseline=(current bounding box.center)]
				\begin{pgfonlayer}{nodelayer}
		\node (0) at (0, 0) {};
		\node (9) at (0, 4.65) {};
	\end{pgfonlayer}
	\begin{pgfonlayer}{edgelayer}
		\draw [style=forward arrow] (9) to (0);
	\end{pgfonlayer}

			\end{tikzpicture}
			%}%
			=
			\,
			%\resizebox{!}{1.0cm}{%
			\begin{tikzpicture}[scale=0.4,baseline=(current bounding box.center)]
			\begin{pgfonlayer}{nodelayer}
		\node [style=none] (0) at (0, 0) {};
		\node [style=none] (1) at (-1, 2) {};
		\node [style=none] (2) at (1, 2) {};
		\node [style=none] (3) at (2, 4) {};
		\node [style=none] (4) at (1, 2) {};
		\node [style=none] (5) at (3, 2) {};
		\node [style=none] (6) at (3, 0) {};
		\node [style=none] (7) at (-1, 4) {};
	\end{pgfonlayer}
	\begin{pgfonlayer}{edgelayer}
		\draw [style=forward arrow] (1.center)
				[in=180, out=-90] to (0.center)
				[in=-90, out=0] to (2.center);
%		\draw [style=forward arrow, in=180, out=-90] (1.center) to (0);
%		\draw [style=forward arrow, in=-90, out=0] (0) to (2.center);
		\draw [style=forward arrow] (4.center)
				[in=-180, out=90] to (3.center)
				[in=90, out=0] to (5.center);
%		\draw [style=forward arrow, in=-180, out=90] (4.center) to (3);
%		\draw [style=forward arrow, in=90, out=0] (3) to (5.center);
		\draw (7.center) to (1.center);
		\draw (5.center) to (6.center);
	\end{pgfonlayer}

			\end{tikzpicture}
			%}%
		\]
\end{Rem}

\begin{Def}\label{def:trace-map}
	Let $A$ be a rigid algebra in the symmetric monoidal category $\cat C$.  Its \emph{trace map} $\tr\colon A\to\unit$ is given by
		\begin{equation}\label{eq:trace-map}
			A \simeq A\otimes \unit \xra{1\otimes \eta} A\otimes DA\otimes A \xra{1\otimes\tau} A\otimes A\otimes DA \xra{\mu\otimes 1} A\otimes DA\xra{\eps}\unit.
		\end{equation}
\end{Def}

\begin{Rem}\label{rem:trace-explanation}
	To explain this definition, recall that every endomorphism $f\colon A \to A$ of the rigid object $A$ has an associated ``trace'' $\Tr(f)\colon \unit \to \unit$ given as
		\[
			\unit \xra{\eta} DA \otimes A \xra{1\otimes f} DA \otimes A \xra{\tau} A \otimes DA \xra{\epsilon} \unit.
		\]
	Moreover, post-composition by the map 
		\begin{equation}\label{eq:blah1}
			DA \otimes A \xra{\tau} A \otimes DA \xra{\epsilon} \unit
		\end{equation}
	induces a function
		\[
			\cat C(A,A) \simeq \cat C(\unit,DA \otimes A) \to \cat C(\unit,\unit)
		\]
	which sends $f$ to $\Tr(f)$. On the other hand, the multiplication map $\mu\colon A \otimes A \to A$ corresponds (by pulling the second $A$ to the other side) to a morphism $A \to DA \otimes A$ given by
		\begin{equation}\label{eq:blah2}
			A \simeq \unit \otimes A \xra{\eta\otimes 1} DA \otimes A \otimes A \xra{1\otimes\tau} DA \otimes A \otimes A \xra{1\otimes \mu} DA \otimes A.
		\end{equation}
	Post-composition by this map provides a function
		\[
			\cat C(\unit, A) \to \cat C(\unit, DA \otimes A) \simeq \cat C(A,A)
		\]
	which sends a morphism $a\colon \unit \to A$ to ``left multiplication by $a$'':
		\[
			L_a \colon  A \simeq A \otimes \unit \xra{1 \otimes a}A \otimes A \xra{\tau} A\otimes A \xra{\mu} A.
		\]
	The map \eqref{eq:trace-map} defining the trace map $\tr\colon A\to \unit$ is readily checked to equal the composite of \eqref{eq:blah2} and~\eqref{eq:blah1}. Post-composition by the trace map thus provides the function
		\begin{align*}
			\cat C(\unit, A)	& \xra{\tr_*} \cat C(\unit,\unit)\\
			a 					& \longmapsto \Tr(L_a).
		\end{align*}
	Thus $\tr\colon A\to\unit$ is morally the map which sends an ``element'' of $A$ to the trace of left multiplication by that element.
\end{Rem}

\begin{Def}\label{def:trace-form}
	The \emph{trace form} of a rigid algebra $A$ is the map $t\colon A \otimes A \to \unit$ defined as the composite
		\[
			A\otimes A \xra{\mu}A \xra{\tr}\unit.
		\]
\end{Def}

\begin{Rem}
	The trace map and trace form of a rigid algebra are given by the following string diagrams:
		\[
			\tr=
			\resizebox{!}{1.5cm}{%
			\begin{tikzpicture}[framed,baseline=(current bounding box.center)]
				\begin{pgfonlayer}{nodelayer}
		\node [style=none] (40) at (-0.75, 3.25) {};
		\node [style=none] (41) at (0.75, 3.25) {};
		\node [style=circle] (44) at (-0.75, 3) {};
		\node [style=none] (49) at (0, 2.5) {};
		\node [style=none] (50) at (0, 2.5) {};
		\node [style=none] (51) at (0.75, 2.5) {};
		\node [style=none] (52) at (-1.5, 2.5) {};
		\node [style=none] (53) at (-1.5, 0.5) {};
		\node [style=none] (54) at (0.75, 1.75) {};
		\node [style=none] (55) at (0, 1.75) {};
		\node [style=none] (56) at (0, 0.75) {};
		\node [style=none] (57) at (0.75, 0.75) {};
	\end{pgfonlayer}
	\begin{pgfonlayer}{edgelayer}
		\draw [in=90, out=0] (44) to (49.center);
		\draw [in=90, out=-90] (50.center) to (54.center);
		\draw [in=-90, out=90, looseness=0.75] (55.center) to (51.center);
		\draw (55.center) to (56.center);
		\draw (54.center) to (57.center);
		\draw (52.center) to (53.center);
		\draw [in=-180, out=90] (52.center) to (44);
		\draw [in=-90, out=-90, looseness=1.25] (56.center) to (57.center);
		\draw (49.center) to (50.center);
		\draw (40.center) to (44);
		\draw (51.center) to (41.center);
		\draw [in=90, out=90, looseness=1.75] (40.center) to (41.center);
	\end{pgfonlayer}

			\end{tikzpicture}}%
			\quad\text{ and }\quad
			t = 
			\resizebox{!}{1.5cm}{%
			\begin{tikzpicture}[framed,baseline=(current bounding box.center)]
				\begin{pgfonlayer}{nodelayer}
		\node [style=none] (40) at (-0.75, 3.25) {};
		\node [style=none] (41) at (0.75, 3.25) {};
		\node [style=circle] (44) at (-0.75, 3) {};
		\node [style=none] (49) at (0, 2.5) {};
		\node [style=none] (50) at (0, 2.5) {};
		\node [style=none] (51) at (0.75, 2.5) {};
		\node [style=none] (52) at (-1.5, 2.5) {};
		\node [style=none] (54) at (0.75, 1.75) {};
		\node [style=none] (55) at (0, 1.75) {};
		\node [style=none] (56) at (0, 1.25) {};
		\node [style=none] (57) at (0.75, 1.25) {};
		\node [style=circle] (58) at (-1.5, 1) {};
		\node [style=none] (59) at (-0.75, 0.5) {};
		\node [style=none] (60) at (-0.75, 0.5) {};
		\node [style=none] (61) at (-2.25, 0.5) {};
		\node [style=none] (62) at (-2.25, -0.5) {};
		\node [style=none] (63) at (-0.75, -0.5) {};
	\end{pgfonlayer}
	\begin{pgfonlayer}{edgelayer}
		\draw [in=90, out=0] (44) to (49.center);
		\draw [in=90, out=-90] (50.center) to (54.center);
		\draw [in=-90, out=90, looseness=0.75] (55.center) to (51.center);
		\draw (55.center) to (56.center);
		\draw (54.center) to (57.center);
		\draw [in=-180, out=90] (52.center) to (44);
		\draw [in=-90, out=-90, looseness=1.25] (56.center) to (57.center);
		\draw (49.center) to (50.center);
		\draw (40.center) to (44);
		\draw (51.center) to (41.center);
		\draw [in=90, out=90, looseness=1.75] (40.center) to (41.center);
		\draw [in=90, out=0] (58) to (59.center);
		\draw [in=-180, out=90] (61.center) to (58);
		\draw (59.center) to (60.center);
		\draw (52.center) to (58);
		\draw (61.center) to (62.center);
		\draw (60.center) to (63.center);
	\end{pgfonlayer}

			\end{tikzpicture}}%
		\]
\end{Rem}

\begin{Rem}\label{rem:invariant-form}
	A map $f\colon A\otimes A \to \unit$ is said to be an ``associative'' form (also called an ``invariant'' form) if $f\circ (\mu \otimes 1) = f\circ (1\otimes \mu)$.  Note that any form which factors through $\mu$ (such as the trace form of a rigid algebra) is necessarily associative by the associativity of the multiplication. The converse is also true: A form $A \otimes A \to \unit$ is associative if and only if it factors through~$\mu$.  In fact, we obtain a bijection
		\[
			\big\{\text{maps } A \to \unit \big\} \isor \big\{ \text{associative forms } A \otimes A \to \unit \big\}
		\]
	given by $\theta \mapsto \theta \circ \mu$
	with inverse $f \mapsto f\circ (u\otimes 1) = f\circ(1\otimes u)$.
\end{Rem}

\begin{Rem}\label{rem:symmetric-form}
	A map $f\colon A \otimes A \to \unit$ is said to be a ``symmetric'' form if $f = f\circ \tau$.  Note that if $A$ is a commutative algebra, then every associative form is automatically symmetric.
\end{Rem}

\begin{Prop}\label{prop:trace-is-symmetric}
	The trace form of a rigid algebra is symmetric.
\end{Prop}

\begin{proof}
	First note that we can rewrite the trace form as
		\begin{equation}\label{eq:string-a}
			\resizebox{!}{\mylengthj}{
			\input{t02.tikz}
			}
			=
			\resizebox{!}{\mylengthj}{
			\input{a02.tikz}
			}
		\end{equation}
	as has already been mentioned in Remark~\ref{rem:trace-explanation}. Next we establish
		\begin{equation}\label{eq:string-lemma}
			\resizebox{!}{\mylengthj}{
			\input{a03.tikz}
			}
			=
			\;\;\;
			\resizebox{!}{\mylengthj}{
			\input{a04.tikz}
			}
		\end{equation}
	which is an equality of morphisms $A \otimes A \to DA \otimes A$.  By adjunction,  this can be checked after applying $A\otimes -$ and post-composing with $A\otimes DA \otimes A \xra{\eps\otimes 1} A$:
		\[
			\resizebox{!}{\mylengthk}{
			\input{a05.tikz}
			}
			\!
			=
			\resizebox{!}{\mylengthk}{
			\input{a06.tikz}
			}
			=
			\resizebox{!}{\mylengthk}{
			\input{a07.tikz}
			}
			=
			\!\!\!
			\resizebox{!}{\mylengthk}{
			\input{a08.tikz}
			}
			=
			\!\!\!\!\!\!
			\resizebox{!}{\mylengthk}{
			\input{a09.tikz}
			}
			=
			\resizebox{!}{\mylengthk}{
			\input{a10.tikz}
			}
			\!\!\!\!
			=
			\resizebox{!}{\mylengthk}{
			\input{a11.tikz}
			}
		\]
	Next note that
		\begin{equation}\label{eq:string-b}
			\resizebox{!}{1cm}{
			\input{a17.tikz}}
			= 
			\resizebox{!}{1cm}{
                        \input{a18.tikz}}
			=
			\;\resizebox{!}{1cm}{
			\input{a19.tikz}
			}
		\end{equation}
	Then
		\[
			\resizebox{!}{\mylengthi}{
			\input{t02.tikz}
			}
			\!\!
			\xeq[\eqref{eq:string-a}]{}
			\resizebox{!}{\mylengthi}{
			\input{a02.tikz}
			}
			\!\!\!\!
			\xeq[\eqref{eq:string-lemma}]{}
			\resizebox{!}{\mylengthi}{
			\input{a12.tikz}
			}
			\!\!\!
			=
			\resizebox{!}{\mylengthi}{
			\input{a13.tikz}
			}
			\xeq[\eqref{eq:string-b}]{}
			\resizebox{!}{\mylengthi}{
			\input{a14.tikz}
			}
			\xeq[\eqref{eq:string-lemma}]{}
			\resizebox{!}{\mylengthi}{
			\input{a15.tikz}
			}
			\!\!
			\xeq[\eqref{eq:string-a}]{}
			\resizebox{!}{\mylengthi}{
			\input{t03.tikz}
			}
		\]
	establishes that the trace form is symmetric.
\end{proof}

\begin{Rem}\label{rem:intuition}
  Intuition for why the trace form is symmetric comes from the fact that for any two endomorphisms $f,g\colon A\to A$, we have $\Tr(f\circ g)=\Tr(g\circ f)$.  Hence, at least morally, $t(a,b) = \Tr(L_{ab}) = \Tr(L_{a} \circ L_{b}) = \Tr(L_{b}\circ L_{a}) = \Tr(L_{ba}) = t(b,a)$.
\end{Rem}

\begin{Def}\label{def:non-degenerate}
	If $A$ is a rigid object in a symmetric monoidal category, then every map $f\colon A\otimes A \to \unit$ gives rise to two morphisms $A \to DA$ by adjunction (moving each copy of $A$ to the right-hand side).  These two maps $A \to DA$ coincide when $f$ is symmetric and are given by
		\begin{equation}\label{eq:f*}
			f^*\colon A\simeq \unit \otimes A \xra{\eta \otimes 1} DA\otimes A \otimes A \xra{1\otimes f} DA \otimes \unit \simeq DA.
		\end{equation}
	We say that a symmetric form $f\colon A\otimes A \to \unit$ is \emph{nondegenerate} if $f^*\colon A \to DA$ is an isomorphism.
\end{Def}

\begin{Prop}\label{prop:strongly-separable-nondegenerate}
	The trace form of a strongly separable rigid algebra is nondegenerate.  Moreover, a strongly separable rigid algebra has a unique symmetric separability idempotent, which is given by
		\begin{equation}\label{eq:kappa-forced}
			\unit \xra{\eta} DA \otimes A \xra{(t^*)^{-1}\otimes 1} A \otimes A.
		\end{equation}
\end{Prop}

\begin{proof}
	Let $\kappa$ be a symmetric separability idempotent. We will start by showing that the composite
		\begin{equation}\label{eq:composite}
			A \simeq \unit \otimes A \xra{\kappa \otimes 1} A \otimes A \otimes A \xra{1 \otimes t} A \otimes \unit \simeq A
		\end{equation}
	is the identity, where $t$ denotes the trace form (Definition~\ref{def:trace-form}). First note:
		\begin{equation}\label{eq:string-c}
			\resizebox{!}{\mylengthd}{
			\input{b10.tikz}
			}
			=
			\resizebox{!}{\mylengthd}{
			\input{b11.tikz}
			}
			=
			\resizebox{!}{\mylengthd}{
			\input{b12.tikz}
			}
		\end{equation}
	Then observe that
		\[
			\resizebox{!}{\mylengthd}{
			\input{t04.tikz}
			}
			\xeq[\kapref2]{}
			%\resizebox{!}{1.5cm}{
			%%\input{b02.tikz}
			%\input{t05.tikz}
			%}
			%=
			\resizebox{!}{\mylengthd}{
			\input{b03.tikz}
			}
			\xeq[\kapref2]{}
			\resizebox{!}{\mylengthd}{
			\input{b04.tikz}
			}
			=
			\resizebox{!}{\mylengthd}{
			\input{b05.tikz}
			}
		\]
	and
		\[
			\resizebox{!}{\mylengthd}{
			\input{b05.tikz}
			}
			=
			\resizebox{!}{\mylengthd}{
			\input{b06.tikz}
			}
			\xeq[\eqref{eq:string-c}]{}
			\resizebox{!}{\mylengthd}{
			\input{b07.tikz}
			}
			\xeq[\kapref4]{}
			\resizebox{!}{\mylengthd}{
			\input{b08.tikz}
			}
			=
			\resizebox{!}{\mylengthd}{
			\begin{tikzpicture}[baseline=(current bounding box.center)]
	\begin{pgfonlayer}{nodelayer}
		\node [style=none] (78) at (2, 2.5) {};
		\node [style=none] (80) at (2, 0.75) {};
	\end{pgfonlayer}
	\begin{pgfonlayer}{edgelayer}
		\draw (78.center) to (80.center);
	\end{pgfonlayer}
\end{tikzpicture}

			}
		\]
	which shows that \eqref{eq:composite} is the identity map.  Now $\kappa$ is symmetric by assumption and the trace form $t$ is symmetric by Proposition~\ref{prop:trace-is-symmetric}. Hence 
		\[
			\resizebox{!}{\mylengthe}{
			\input{c01.tikz}
			}
			=
			\resizebox{!}{\mylengthe}{
			\input{c02.tikz}
			}
			=
			\resizebox{!}{\mylengthe}{
			\input{c03.tikz}
			}
			=
			\resizebox{!}{\mylengthe}{
			\input{c04.tikz}
			}
		\]
	In other words, the composite \eqref{eq:composite} coincides with the other composite
		\[
			A \simeq A \otimes \unit \xra{1\otimes \kappa} A \otimes A \otimes A \xra{t \otimes 1} \unit \otimes A \simeq A.
		\]
	It follows that the map $\kappa^*\colon DA \to A$ given by 
		\[
			DA \simeq \unit \otimes DA \xra{\kappa \otimes 1} A \otimes A \otimes DA \xra{1 \otimes \epsilon} A \otimes \unit \simeq A
		\]
	is an inverse to $t^*\colon A \to DA$. Indeed: 
		\[ 
			\resizebox{!}{\mylengthe}{
			\input{d01.tikz}
			}
			=
			\;
			\resizebox{!}{\mylengthe}{
			\input{d02.tikz}
			}
			%\;
			=
			\;
			\resizebox{!}{\mylengthe}
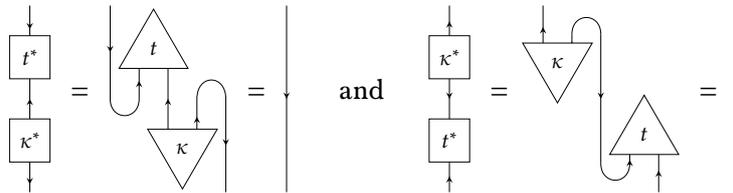
{
			\begin{tikzpicture}[scale=0.5,baseline=(current bounding box.center)]
	\begin{pgfonlayer}{nodelayer}
		\node [style=none] (0) at (0, 3.5) {};
		\node [style=none] (1) at (0, -3) {};
	\end{pgfonlayer}
	\begin{pgfonlayer}{edgelayer}
		\draw [style=forward arrow] (0.center) to (1.center);
	\end{pgfonlayer}
\end{tikzpicture}

			}
			\quad
			\text{ and }
			\quad
			\resizebox{!}{\mylengthe}{
			\input{d04.tikz}
			}
			=
			%\;
			\resizebox{!}{\mylengthe}{
			\input{d05.tikz}
			}
			%\;
			=
			\;
			\resizebox{!}{\mylengthe}{
			\begin{tikzpicture}[scale=0.5,baseline=(current bounding box.center)]
	\begin{pgfonlayer}{nodelayer}
		\node [style=none] (0) at (0, 3.5) {};
		\node [style=none] (1) at (0, -3) {};
	\end{pgfonlayer}
	\begin{pgfonlayer}{edgelayer}
		\draw [style=backward arrow] (0.center) to (1.center);
	\end{pgfonlayer}
\end{tikzpicture}

			}
		\]
	In particular, the trace form is nondegenerate.  Moreover, one can readily check that $(t^*\otimes 1)\circ \kappa = \eta$, from which it follows that $\kappa$ is given by \eqref{eq:kappa-forced}. 
\end{proof}

\begin{Thm}\label{thm:strongly-separable-characterization}
	A rigid algebra is strongly separable if and only if its trace form is nondegenerate.
\end{Thm}

\begin{proof}
	The ``only if'' part is provided by Proposition~\ref{prop:strongly-separable-nondegenerate}. Conversely, suppose $A$ is a rigid algebra whose trace form $t\colon A \otimes A \to \unit$ is nondegenerate. Write $\theta\colon A \isor DA$ for the associated isomorphism (that is, $\theta =t^*$ in the notation of Definition~\ref{def:non-degenerate}), and define $\kappa\colon \unit \to A \otimes A$ by 
		\begin{equation}\label{eq:kdef}
			\unit \xra{\eta} DA \otimes A \xra{\theta^{-1}\otimes 1} A \otimes A.
		\end{equation}
	First we check that $t$ and $\kappa$ form a self-duality in the sense that
		\begin{equation}\label{eq:self-dual-1}
			\resizebox{!}{\mylengthh}{
			\input{f04.tikz}
			}
			\xeq[\eqref{eq:kdef}]{}
			\resizebox{!}{\mylengthh}{
			\input{f03.tikz}
			}
			\xeq[\eqref{eq:f*}]{}
			\;
			\resizebox{!}{\mylengthh}{
			\input{f01.tikz}
			}
			=
			\;
			\resizebox{!}{\mylengthh}{
			\begin{tikzpicture}[scale=0.7, baseline={(current bounding box.center)}]
	\begin{pgfonlayer}{nodelayer}
		\node [style=none] (0) at (0.5, -3.75) {};
		\node [style=none] (6) at (0.5, 3.75) {};
	\end{pgfonlayer}
	\begin{pgfonlayer}{edgelayer}
		\draw [style=backward arrow] (6.center) to (0.center);
	\end{pgfonlayer}
\end{tikzpicture}

			}
		\end{equation}
	and
		\begin{equation}\label{eq:self-dual-2}
			\resizebox{!}{\mylengthh}{
			\input{f05.tikz}
			}
			%\;	
			\xeq[\eqref{eq:kdef}]{}
			\;
			\resizebox{!}{\mylengthh}{
			\input{f07.tikz}
			}
			\;	
			\xeq[\eqref{eq:f*}]{}
			\;
			%\;\;
			\resizebox{!}{\mylengthh}{
			\input{f06.tikz}
			}
			\;
			=
			\;
			\resizebox{!}{\mylengthh}{
			\begin{tikzpicture}[scale=0.7, baseline={(current bounding box.center)}]
	\begin{pgfonlayer}{nodelayer}
		\node [style=none] (0) at (0.5, -3.75) {};
		\node [style=none] (6) at (0.5, 3.75) {};
	\end{pgfonlayer}
	\begin{pgfonlayer}{edgelayer}
		\draw [style=backward arrow] (6.center) to (0.center);
	\end{pgfonlayer}
\end{tikzpicture}

			}
			%=
			%\;
			%\input{f02.tikz}
		\end{equation}
	Armed with this relationship between $t$ and $\kappa$, the fact that $t$ is symmetric (by Proposition~\ref{prop:trace-is-symmetric}) implies that $\kappa$ is also symmetric:
		\[
			\resizebox{\textwidth}{!}{
			\input{e01.tikz}
			$\xeq[\eqref{eq:self-dual-1}]{}$
			\input{e02.tikz}
			$\xeq[(\ref{prop:trace-is-symmetric})]{}$
			\input{e03.tikz}
			=
			\input{e04.tikz}
			$\xeq[\eqref{eq:self-dual-1}]{}$
			\input{e05.tikz}
			=
			\input{e06.tikz}
			}
		\]
	This establishes axiom \kapref3.  Next we establish \kapref2, visualized in string diagrams in \eqref{eq:string-k2}. It suffices to check equality after post-composition by the isomorphism $\theta \otimes \id_A$. Then by adjunction, it suffices to check equality after applying $A\otimes -$ and post-composing by $A\otimes DA \otimes A \xra{\epsilon \otimes 1} A$. Indeed, 
		\[
			\resizebox{\textwidth}{!}{
			\input{g02.tikz}
			%$\underset{\text{def}}{=}\,$
			$\!\xeq[\eqref{eq:kdef}]{}\,$
			\input{g01.tikz}
			=
			\input{g03.tikz}
			$\xeq[\eqref{eq:self-dual-2}]{}$
			\input{g04.tikz}
			$\xeq[(\dagger)]{}$
			\input{g05.tikz}
			$
			\!\!\!\!\!
			%\underset{\text{def}}{=}\,$
			\xeq[\eqref{eq:f*}]{}\,$
			\input{g06.tikz}
			}
		\]
	where the equality $(\dagger)$ is the fact that the trace form is an 
	associative
	form (Remark~\ref{rem:invariant-form}). Finally, we establish~\kapref1. Observe that
		\[
			\resizebox{!}{\mylengthg}{
			\begin{tikzpicture}[baseline=(current bounding box.center)]
	\begin{pgfonlayer}{nodelayer}
		\node [style=none] (0) at (4.25, 4) {};
		\node [style=none] (1) at (4.25, -1.5) {};
	\end{pgfonlayer}
	\begin{pgfonlayer}{edgelayer}
		\draw [in=270, out=90] (1.center) to (0.center);
	\end{pgfonlayer}
\end{tikzpicture}

			}
			=
			\resizebox{!}{\mylengthg}{
			\input{h03.tikz}
			}
			=
			\resizebox{!}{\mylengthg}{
			\input{h08.tikz}
			}
			\xeq[(\ddagger)]{}
			\resizebox{!}{\mylengthg}{
			\input{h09.tikz}
			}
			=
			\resizebox{!}{\mylengthg}{
			\input{h01.tikz}
			}
			=
			\resizebox{!}{\mylengthg}{
			\input{h02.tikz}
			}
		\]
	and note that we showed $(\ddagger)$ was a consequence of \kapref2 in the proof of Proposition~\ref{prop:strongly-separable-nondegenerate}. Precomposing with the unit, we obtain
		\[
			\resizebox{!}{\mylengthg}{
			\begin{tikzpicture}[baseline=(current bounding box.center)]
	\begin{pgfonlayer}{nodelayer}
		\node [style=none] (0) at (4.25, 4) {};
		\node [style=tiny circle] (1) at (4.25, -1.5) {};
	\end{pgfonlayer}
	\begin{pgfonlayer}{edgelayer}
		\draw [in=270, out=90] (1) to (0.center);
	\end{pgfonlayer}
\end{tikzpicture}

			}
			=
			\resizebox{!}{\mylengthg}{
			\input{h06.tikz}
			}
			=
			\resizebox{!}{\mylengthg}{
			\input{h07.tikz}
			}
		\]
	which is \kapref1.
\end{proof}

\begin{Cor}\label{cor:rigid-comm-alg}
	A rigid commutative algebra is separable if and only if it is strongly separable if and only if its trace form is nondegenerate.
\end{Cor}

\begin{proof}
	Every commutative separable algebra is strongly separable (Corollary~\ref{cor:commutative-separable-is-strongly-separable}); hence the claim follows from Theorem~\ref{thm:strongly-separable-characterization}.
\end{proof}

\begin{Rem}\label{rem:rigid-strongly-separable-is-self-dual}
	A rigid strongly separable algebra $A$ is automatically self-dual since the nondegeneracy of the trace form provides an isomorphism $A \cong DA$.
\end{Rem}

\begin{Exa}
	Consider the case where $\cat C=R\MMod$ is the category of $R$-modules for $R$ a commutative ring. An $R$-algebra $A$ is rigid precisely when it is finitely generated and projective (equivalently, finitely presented and flat) as an $R$-module. The trace map $A\to R$ is $a \mapsto \Tr(L_a)$,  where $L_a \colon  A \to A$ denotes left multiplication by $a$, and the trace form $t\colon A\otimes A \to R$ is given by $t(a\otimes b) = \Tr(L_{ab})$. In this example, the argument in Remark~\ref{rem:intuition} shows immediately that the trace form is symmetric.  It turns out that over a field $R=k$, a separable algebra is automatically rigid (that is, finite-dimensional), as shown by \cite[Proposition~1.1]{VillamayorZelinsky66}.  It was partly to clarify such finiteness assumptions that led the author to write this section on strongly separable algebras in arbitrary symmetric monoidal categories.
\end{Exa}

\begin{Exa}
	An \emph{idempotent algebra} in a symmetric monoidal category is an algebra $(A,\mu,u)$ whose multiplication map $\mu\colon A \otimes A \to A$ is an isomorphism. This is equivalent to the equality $u\otimes A=A \otimes u$ of morphisms $A \to A\otimes A$ (which then serve as an inverse to $\mu$). It is also equivalent to the switch map $\tau \colon A \otimes A \to A \otimes A$ being equal to the identity map $A \otimes A \to A \otimes A$. Idempotent algebras are thus examples of commutative (strongly) separable algebras. They have a (unique) separability idempotent given by $\mu^{-1} \circ u\colon \unit \to A \otimes A$. However, they are usually not rigid. For example, take $\cat C=R\MMod$ for $R$ a commutative ring. The idempotent \mbox{$R$-algebra}~$R[1/s]$ is rarely finitely generated as an $R$-module. Indeed, this would imply that the principal open $D(s) \subset \Spec(R)$ is both an open and closed subset of $\Spec(R)$; see the argument in \cite[Example~7.4]{Sanders19}, for example.
\end{Exa}

\begin{Rem}
	A discussion of separable algebras would not be complete without saying something about their relationship with Frobenius algebras: 
\end{Rem}

\begin{Def}
	A \emph{Frobenius algebra} in a symmetric monoidal category is an object~$A$ equipped with both an algebra structure $(A,\mu,u)$ and a coalgebra structure $(A,\Delta,c)$ such that the Frobenius law holds:
		\[
			(1\otimes \mu)\circ(\Delta\otimes 1) = \Delta \circ \mu = (\mu \otimes 1)\circ (1 \otimes \Delta).
		\]
	See, for example, \cite[3.6.8]{Kock04}. We say that $A$ is a \emph{symmetric} Frobenius algebra if the associative form
		\[
			A \otimes A \xra{\mu} A \xra{c} \unit
		\]
	is symmetric. Thus, every commutative Frobenius algebra is symmetric. A \emph{special} Frobenius algebra is a Frobenius algebra such that $\mu \circ \Delta = \id_A$.
\end{Def}

\begin{Rem}\label{rem:frobenius-self-dual}
	If $(A,\mu,u,\Delta,c)$ is a Frobenius algebra, then the underlying object~$A$ is necessarily self-dual (\textit{cf.}~Remark~\ref{rem:rigid-strongly-separable-is-self-dual}).  Indeed, the two maps $c\circ \mu \colon  A \otimes A \to \unit$ and $\Delta\circ u\colon \unit \to A\otimes A$ provide a self-duality.
\end{Rem}

\begin{Rem}
	The following relationship between strongly separable algebras and special symmetric Frobenius algebras is well-known classically; we include a proof for precision and completeness.  The interested reader will find more concerning these ideas in \cite[Section 2.5]{LaudaPfeiffer07}, \cite[Section 3.3]{FuchsRunkelSchweigert02}, and \cite{Fauser13}, among other sources.
\end{Rem}

\begin{Prop}
	An algebra admits the structure of a special symmetric Frobenius algebra if and only if it is rigid and strongly separable.  In this case, the special symmetric Frobenius structure is unique: The counit $A \to \unit$ is the trace map $($Definition~\ref{def:trace-map}\,$)$, and the comultiplication $A \to A \otimes A$ is the map corresponding $($Remark~\ref{rem:precompose}\,$)$ to the unique symmetric separability idempotent $($Proposition~\ref{prop:strongly-separable-nondegenerate}\,$)$.
\end{Prop}

\begin{proof}
	If $A$ is a strongly separable rigid algebra with (unique) symmetric separability idempotent $\kappa\colon \unit \to A \otimes A$, then the corresponding map $A \to A \otimes A$ is coassociative.  Indeed, using both descriptions provided by \kapref{2}, we have  
		\[
			\resizebox{!}{\mylength}{
			\input{s01.tikz}
			}
			%\underset{\eqref{eq:string-a}}{=}
			%=
			\xeq[\kapref2]{}
			\resizebox{!}{\mylength}{
			\input{s02.tikz}
			}
			%\underset{\eqref{eq:string-lemma}}{=}
			\xeq[\kapref2]{}
			%=
			\resizebox{!}{\mylength}{
			\input{s03.tikz}
			}
			%\underset{(\text{assoc})}{=}
			=
			\resizebox{!}{\mylength}{
			\input{s04.tikz}
			}
			%=
			\xeq[\kapref2]{}
			\resizebox{!}{\mylength}{
			\input{s05.tikz}
			}
			%=
			\xeq[\kapref2]{}
			\resizebox{!}{\mylength}{
			\input{s06.tikz}
			}
		\]
	This provides $A$ with the structure of a coalgebra with counit $A \to \unit$ given by the trace map.  For the counital axiom, just observe that
		\[
			\resizebox{!}{\mylengthb}{
			\input{u01.tikz}
			}
			%\underset{\eqref{eq:string-a}}{=}
			%=
			%\resizebox{!}{\mylengthb}{
			%\input{u02.tikz}
			%}
			%\underset{\eqref{eq:string-lemma}}{=}
			\xeq[\kapref2]{}
			\resizebox{!}{\mylengthb}{
			\input{u03.tikz}
			}
			=
			\resizebox{!}{\mylengthb}{
			\input{u04.tikz}
			}
			%=
			\xeq[(\dagger)]{}
			\resizebox{!}{\mylengthb}{
			\begin{tikzpicture}
	\begin{pgfonlayer}{nodelayer}
		\node [style=none] (3) at (0.5, 6) {};
		\node [style=none] (4) at (0.5, -3) {};
	\end{pgfonlayer}
	\begin{pgfonlayer}{edgelayer}
		\draw (3.center) to (4.center);
	\end{pgfonlayer}
\end{tikzpicture}

			}
			\quad\text{ and }\quad
			\resizebox{!}{\mylengthb}{
			\input{u06.tikz}
			}
			%=
			\xeq[\kapref2]{}
			\resizebox{!}{\mylengthb}{
			\input{u07.tikz}
			}
			=
			\resizebox{!}{\mylengthb}{
			\input{u08.tikz}
			}
			%=
			\xeq[(\dagger)]{}
			\resizebox{!}{\mylengthb}{
			\begin{tikzpicture}
	\begin{pgfonlayer}{nodelayer}
		\node [style=none] (3) at (0.5, 6) {};
		\node [style=none] (4) at (0.5, -3) {};
	\end{pgfonlayer}
	\begin{pgfonlayer}{edgelayer}
		\draw (3.center) to (4.center);
	\end{pgfonlayer}
\end{tikzpicture}

			}
		\]
	where the last equalities $(\dagger)$ were established in the proof of Proposition~\ref{prop:strongly-separable-nondegenerate}.  Alternatively, one can use the description~\eqref{eq:kappa-forced} of the unique separability idempotent and check the counital diagrams after post-composition by the isomorphism $t^*\colon A \isor DA$.  This establishes that a strongly separable rigid algebra admits the structure of a special symmetric Frobenius algebra.

	Now suppose that $(A,\mu,u,\Delta,c)$ is a special symmetric Frobenius algebra.  Every Frobenius algebra is self-dual (Remark~\ref{rem:frobenius-self-dual}), and the comultiplication $\Delta\colon A \to A \otimes A$ satisfies \sigref{2}.  In our case, it also satisfies \sigref{1} since $A$ is assumed to be special. Symmetry of the associated separability idempotent $\Delta\circ u\colon  \unit \to A \otimes A$ then follows from the assumed symmetry of $c\circ \mu\colon A\otimes A \to \unit$ via the self-duality (as in the beginning of the proof of Theorem~\ref{thm:strongly-separable-characterization}).  Thus $A$ is strongly separable with symmetric separability idempotent $\Delta \circ u$.

	To establish uniqueness, first observe that the commutative diagram
		\[\begin{tikzcd}[column sep=large]
				A  \ar[d,"\Delta"] \ar[r,"u \otimes 1"] & A\otimes A \ar[d,"1\otimes \Delta"] \ar[r,"\Delta \otimes 1"] & A\otimes A \otimes A \ar[d,"1\otimes \mu"] \\
				A \otimes A \ar[rr,bend right=15,"\id"'] \ar[r,"u\otimes 1\otimes 1"] & A \otimes A \otimes A \ar[r,"\mu \otimes 1"] & A \otimes A
		\end{tikzcd}\]
	shows that the comultiplication $\Delta$ of a Frobenius algebra is determined by $\Delta \circ u$ and $\mu$.  If $(A,\mu,u,\Delta_1,c_1)$ and $(A,\mu,u,\Delta_2,c_2)$ are two special symmetric Frobenius structures on the (rigid strongly separable) algebra $(A,\mu,u)$,  then $\Delta_1 \circ u = \Delta_2 \circ u$ by the uniqueness of symmetric separability idempotents (Proposition~\ref{prop:strongly-separable-nondegenerate}), and hence $\Delta_1 = \Delta_2$.  Moreover, since $c_i\circ \mu\colon  A \otimes A \to \unit$ and $\Delta_i \circ u\colon \unit \to A \otimes A$ form a self-duality for each $i=1,2$, we have: 
		\[
			\resizebox{!}{\mylengthc}{
			\input{u09.tikz}
			}
			=
			\resizebox{!}{\mylengthc}{
			\input{u10.tikz}
			}
			=
			\resizebox{!}{\mylengthc}{
			\input{u11.tikz}.
			}
		\]
	That is, $c_1 \circ \mu = c_2 \circ \mu$.  Precomposing by the unit, we conclude that $c_1=c_2$.  This establishes that an algebra admits at most one special symmetric Frobenius structure.  Finally, we have already proved that if $A$ admits a special symmetric Frobenius structure then it is rigid and strongly separable and consequently the trace map and the symmetric separability idempotent provide it with a special symmetric Frobenius structure. These thus provide the unique such structure.
\end{proof}

\section{Separable algebras and triangulated categories}\label{sec:separable-triangulated}

In this section, we recall the relationship between separable algebras and tensor-triangulated categories established in \cite{Balmer11}.

\begin{Rem}
	Recall from \cite[Section 5]{Balmer11} that for each $2\le N \le \infty$, there is the notion of an $N$-triangulated category (or triangulated category of order $N$) which includes as part of the structure a distinguished class of $n$-triangles for each ${n \le N}$ which are required to satisfy suitable higher octahedral axioms. A $2$-triangulated category is precisely the same thing as a pre-triangulated category, while the usual notion of triangulated category (in the sense of Verdier) lies between the notion of \mbox{$2$-triangulated} and $3$-triangulated. An $N$-triangulated functor is a functor which commutes with the suspension and preserves distinguished $N$-triangles (equivalently, preserves distinguished $n$-triangles for all $n \le N$).
\end{Rem}

\begin{Exa}
	The homotopy category $\Ho(\cat C)$ of a stable $\infty$-category has the structure of an $\infty$-triangulated category.
\end{Exa}

\begin{Rem}
	A tensor-triangulated category is a triangulated category equipped with a closed symmetric monoidal structure which is compatible with the triangulation in the sense of \cite[Definition A.2.1]{HoveyPalmieriStrickland97}. For $2\le N \le \infty$, we similarly have the notion of an $N$-tensor-triangulated category by replacing all instances of ``triangulated'' in the definition with ``$N$-triangulated.'' By an ($N$-)tensor-triangulated functor, we mean an ($N$-)triangulated functor which is also a strong symmetric monoidal functor.
\end{Rem}

\begin{Exa}
	The homotopy category $\Ho(\cat C)$ of a presentably symmetric monoidal \cite[Definition~2.1]{NikolausSagave17} stable $\infty$-category is an $\infty$-tensor-triangulated category.
\end{Exa}

\begin{Exa}\label{exa:A-modules-triangulated}
	If $A$ is a commutative separable algebra in an $N$-tensor-triangulated category $\cat T$ ($2 \le N \le \infty$), then the Eilenberg--Moore category $A\MMod_{\cat T}$ inherits the structure of an $N$-tensor-triangulated category such that the extension-of-scalars functor $F_A \colon  \cat T \to A\MMod_{\cat T}$ is an $N$-tensor-triangulated functor. The distinguished $n$-triangles in $A\MMod_{\cat T}$ ($n \le N$) are precisely those which are created by the forgetful functor $U_A \colon A\MMod_{\cat T} \to \cat T$. This is established by \cite[Theorem~5.17]{Balmer11} and \cite[Section 1]{Balmer14}.
\end{Exa}

\begin{Rem}\label{rem:DS-monadic}
	The main theorem of \cite{DellAmbrogioSanders18} states that if $\cat T$ is an idempotent complete triangulated category, then any triangulated adjunction $F\colon \cat T \adjto \cat S:G$ is essentially monadic (that is, monadic up to idempotent completion and killing the kernel of~$G$) whenever the Eilenberg--Moore category inherits a triangulation from~$\cat T$:
		\begin{equation}\label{eq:monadic-equiv}
			(\cat S/\ker G)^\natural \cong GF\MMod_{\cat T}.
		\end{equation}
	This theorem also holds (with the same proof) in the 2-category of $N$-triangulated categories for any $2 \le N \le \infty$. In this case, the equivalence~\eqref{eq:monadic-equiv} is an equivalence of $N$-triangulated categories. To be clear, this is under the hypothesis that the Eilenberg--Moore category $GF\MMod_{\cat T}$ inherits an $N$-triangulation from the \mbox{$N$-triangulation} of $\cat T$ (see \cite[Remark 1.8]{DellAmbrogioSanders18}). This is a strong hypothesis on the adjunction but, as established by Balmer (Example~\ref{exa:A-modules-triangulated}), does hold in the separable case. The following proposition clarifies the situation with the tensor:
\end{Rem}

\begin{Prop}\label{prop:tensor-monadic}
  Let $F\colon\cat T \to \cat S$ be an $(N$-$)$tensor-triangulated functor with $\cat T$ idempotent complete. Suppose $F$ admits a right adjoint $G$ whose kernel is closed under the tensor product, and that the $F \dashv G$ adjunction satisfies the right projection formula \cite[Definition~2.7]{BalmerDellAmbrogioSanders15}.
  If the commutative algebra ${G(\unit) \in \cat T}$ is separable then we have an induced equivalence
		\[
			(\cat S/\ker G)^\natural \cong G(\unit)\MMod_{\cat T}
		\]
	of $(N$-$)$tensor-triangulated categories.
\end{Prop}

\begin{proof}
	By \cite[Lemma~2.8]{BalmerDellAmbrogioSanders15}, the projection formula implies that the monad of the adjunction is the monad associated to the algebra $G(\unit)$. Since $G(\unit)$ is separable, the Eilenberg--Moore category inherits a triangulation from $\cat T$ (Example~\ref{exa:A-modules-triangulated}); hence by \cite{DellAmbrogioSanders18} and Remark~\ref{rem:DS-monadic}, we have equivalences
		\begin{equation}\label{eq:separable-monadic-equivs}
			(G(\unit)\FFree_{\cat T})^\natural \isor (\cat S/\ker G)^\natural \isor G(\unit)\MMod_{\cat T}.
		\end{equation}

	We next prove that $\ker G$ is a tensor-ideal. To this end, let $q\colon\cat S \to \cat S/\ker G$ denote the Verdier quotient and $\overbar{G}\colon \cat S/\ker G \to \cat T$ the induced functor. The induced adjunction $q\circ F \dashv \overbar{G}$ realizes the same monad as the $F \dashv G$ adjunction. Since the comparison functor $K\colon\cat S/\ker G \to G(\unit)\MMod_{\cat T}$ is fully faithful (hence separable), the functor $\overbar{G} = U\circ K$ is separable, being the composite of such (\textit{cf.}~\cite[Proposition~3.5]{Chen15}). Hence, the counit of the $q\circ F \dashv \overbar{G}$ adjunction splits. Thus, for any $s \in \cat S$, the morphism $\alpha$ in the exact triangle $FGs \xra{\eps_s} s \xra{\alpha} x \to \Sigma FGs$ satisfies $q(\alpha)=0$. This means that $\beta \circ \alpha = 0$ for some morphism $\beta\colon  x \to y$ with $\cone(\beta) \in \ker G$. It follows that $s \in \thick\langle FGs,\ker G\rangle$. We conclude that $\ker G$ is a tensor-ideal since the thick subcategory $\{a \in \cat S\mid a\otimes \ker G \subseteq \ker G\}$ contains $\ker G$ by the hypothesis that $\ker G$ is closed under the tensor product, and it contains the image of $F$ by the projection formula.

	We have established that the thick subcategory $\ker G$ is a tensor-ideal. Thus $\cat S/\ker G$ and its idempotent completion $(\cat S/\ker G)^\natural$ inherit tensor structures from~$\cat S$. On the other hand, the Kleisli category $G(\unit)\FFree_{\cat T}$ inherits a tensor structure from~$\cat T$ such that the canonical functor $\cat T \to G(\unit)\FFree_{\cat T}$ is a strict symmetric monoidal functor. The canonical functor $G(\unit)\FFree_{\cat T} \to \cat S$ then inherits the structure of a strong symmetric monoidal functor from the corresponding structure on~$F$. The first functor in \eqref{eq:separable-monadic-equivs} is a strong symmetric monoidal equivalence. It follows that the second functor is also a symmetric monoidal equivalence since the tensor structure on $G(\unit)\MMod_{\cat T}\cong(G(\unit)\FFree_{\cat T})^\natural$ is the idempotent completion of the tensor structure on the Kleisli category (see \cite[Section~1.1]{Pauwels_dissertation} and \cite[Section 1]{Balmer14}).
\end{proof}

\begin{Rem}\label{rem:thm-1}
	As an application of the proposition, we can provide a proof of Theorem~\ref{thm:smashing} stated in the Introduction which characterizes smashing localizations of rigidly-compactly generated categories.
\end{Rem}

\begin{proof}[Proof of Theorem~\ref{thm:smashing}]
  The $(\Rightarrow)$ direction is well-known: Any smashing localization of a rigidly-compactly generated category is a geometric functor to a rigidly-compactly generated category (whose right adjoint is fully faithful); see \cite[Section~3.3]{HoveyPalmieriStrickland97}. For the $(\Leftarrow)$ direction, recall that smashing localizations are nothing but extension-of-scalars with respect to idempotent algebras. Suppose $f^*\colon\cat D \to \cat C$ is a geometric functor whose right adjoint $f_*$ is fully faithful. The multiplication map $f_*(\unitC)\otimes f_*(\unitC) \to f_*(\unitC)$ becomes, under the projection formula, the counit $f_*(\eps)\colon  f_*f^*f_*(\unitC) \simeq f_*(\unitD \otimes f^*f_*(\unitC)) \simeq f_*(\unitC) \otimes f_*(\unitC) \to f_*(\unitC)$. This is an isomorphism since $f_*$ is fully faithful. Thus $f_*(\unitC)$ is an idempotent algebra. Idempotent algebras are separable, so Proposition~\ref{prop:tensor-monadic} gives the result.
\end{proof}

\begin{Rem}
	Note that we assume in Proposition~\ref{prop:tensor-monadic} that the kernel of the right adjoint $G$ is closed under the tensor product and show, under the other hypotheses, that it is then necessarily a tensor-ideal. The next example clarifies that this hypothesis on the kernel of $G$ is not forced by the other assumptions.
\end{Rem}

\begin{Exa}\label{exa:loc}
	Let $\cat T$ be a rigidly-compactly generated tensor-triangulated category. Then $\cat L\coloneqq \Loc(\unit)$ is a tensor-triangulated subcategory of $\cat T$. It is also rigidly-compactly generated, and the inclusion functor $F\colon \cat L \hookrightarrow \cat T$ is a coproduct-preserving tensor-triangulated functor. Hence it has a right adjoint $G$ and the projection formula holds. Moreover, since the left adjoint $F$ is fully faithful, $G(\unit) \simeq \unit$ is just the trivial ring, which is certainly separable. Since the inclusion $\cat L \hookrightarrow \cat T$ has a right adjoint, $\cat L$ is the kernel of a Bousfield localization on $\cat T$. The image of this Bousfield localization is $\cat L^\perp = \ker G$. Hence $\cat L = {}^\perp(\cat L^\perp) = {}^\perp(\ker G)$. If $\ker G$ were a tensor-ideal, then $\cat L$ would be forced to be a tensor-ideal. Indeed, the localizing subcategory $\cat L \subseteq \cat T$ is a tensor-ideal if and only if it is closed under tensoring with compact(=rigid) objects, and the left orthogonal of a tensor-ideal is certainly closed under tensoring with rigid objects. Thus, if $\ker G$ is a tensor-ideal, then $\Loc(\unit)$ is a tensor-ideal, and this is the case if and only if $\Loc(\unit) = \cat T$. Thus, if we take $\cat T$ to be any non-monogenic rigidly-compactly generated tensor-triangulated category, we conclude that $\ker G$ is not a tensor-ideal, and in fact not closed under the tensor product.
\end{Exa}

\section{Finite \'{e}tale morphisms}\label{sec:finite-etale}

The idea that extension-of-scalars with respect to a commutative separable algebra provides tensor triangular geometry with an analogue of an \'{e}tale extension goes back to the work of Balmer \cite{Balmer15,Balmer16a,Balmer16b}. Here we focus on \emph{finite} \'{e}tale extensions of rigidly-compactly generated categories.

\begin{Ter}
	A coproduct-preserving ($N$-)tensor-triangulated functor between rigidly-compactly generated ($N$-)tensor-triangulated categories is called a \emph{geometric functor}.
\end{Ter}

\begin{Def}\label{def:finite-etale}
	A geometric functor $f^*\colon \cat D \to \cat C$ between rigidly-compactly generated ($N$-)tensor-triangulated categories is \emph{finite \'{e}tale} if there exists a compact commutative separable algebra $A$ in $\cat D$ and an ($N$-)tensor-triangulated equivalence\footnote{\, For tensor-triangulated categories in the usual sense of Verdier, the category of modules $A\MMod_{\cat D}$ is \textit{a priori} only a pre-tensor-triangulated category, but this does not cause any trouble for the definition.  Since $\cat C$ is tensor-triangulated by assumption, the equivalence $\cat C \cong A\MMod_{\cat D}$ just forces $A\MMod_{\cat D}$ to be tensor-triangulated as well.  This technicality doesn't arise when working in the 2-category of $N$-tensor-triangulated categories for any $2 \le N \le \infty$.} $\cat C \cong A\MMod_{\cat D}$ such that the functor $f^*$ becomes isomorphic to the extension-of-scalars functor $F_A \colon \cat D \to A\MMod_{\cat D}$.
\end{Def}

\begin{Rem}\label{rem:AModisrcg}
	It follows from \cite[Theorem~4.2]{Balmer16a} that if $\cat D$ is rigidly-compactly generated, then $A\MMod_{\cat D}$ is also rigidly-compactly generated (for $A$ a commutative separable algebra in $\cat D$). Thus, there is no loss of generality in considering only geometric functors between rigidly-compactly generated categories.
\end{Rem}

\begin{Rem}
	Recall from \cite{BalmerDellAmbrogioSanders16} that a geometric functor $f^*\colon\cat D \to \cat C$ between rigidly-compactly generated tensor-triangulated categories has a right adjoint $f_*\colon\cat C \to \cat D$ which itself has a right adjoint $f^!\colon\cat D \to \cat C$. The \emph{relative dualizing object} of~$f^*$ is the object $\omega_f \coloneqq f^!(\unit_{\cat D}) \in \cat C$. Recall that $f^*$ is said to satisfy \emph{Grothendieck--Neeman duality} if the right adjoint $f_*$ preserves compact objects. (A number of equivalent definitions are provided by \cite[Theorem~3.3]{BalmerDellAmbrogioSanders16}.) In this case, the commutative algebra $f_*(\unitC)$ is compact(=rigid). Hence it has a trace map ${f_*(\unitC) \to \unitD}$ (Definition~\ref{def:trace-map}), which corresponds to a map $\unitC \to \omega_f$.
\end{Rem}

\begin{Rem}
	In general, morphisms $\unitC \to \omega_f$ can be identified with morphisms $f_*(\unitC) \to \unitD$ by adjunction, and these can be identified as in Remark~\ref{rem:invariant-form} with the associative forms on the algebra $f_*(\unitC)$:
	\[
		\big\{ \unitC \to \omega_f \big\} \isor
		\big\{ f_*(\unitC) \to \unitD \big\}
		\isor 
		\big\{\text{associative forms }
		f_*(\unitC)\otimes f_*(\unitC) \to \unitD\big\}.
	\]
	Also recall (Definition~\ref{def:non-degenerate}) that an associative form $f_*(\unitC) \otimes f_*(\unitC) \to \unitD$ is nondegenerate if the adjoint morphism $f_*(\unitC) \to Df_*(\unitC)$ is an isomorphism. (Note that these associative forms are automatically symmetric since the algebra $f_*(\unitC)$ is commutative.) On the other hand, recall from \cite[(2.18)]{BalmerDellAmbrogioSanders16} that we have an isomorphism $f_*(\omega_f) \simeq Df_*(\unitC)$.
\end{Rem}

\begin{Lem}\label{lem:omega-form}
	Let $\theta\colon\unitC \to \omega_f$ be any morphism. The map $f_*(\unitC)\to Df_*(\unitC)$ which is adjoint to the associative form on $f_*(\unitC)$ corresponding to $\theta$ coincides with the map
		\begin{equation}\label{eq:omega-dual}
			f_*(\unitC) \xra{f_*(\theta)} f_*(\omega_f) \simeq Df_*(\unitC).
		\end{equation}
	Consequently, the associative form associated to $\theta$ is nondegenerate if and only if $f_*(\theta)$ is an isomorphism.
\end{Lem}

\begin{proof}
  This is a straightforward verification. From the definition of the isomorphism $f_*(\omega_f)\simeq Df_*(\unitC)=[f_*(\unitC),\unitD]$ in \cite[(2.18)]{BalmerDellAmbrogioSanders16}, one sees that the morphism~\eqref{eq:omega-dual} is obtained by going along the top of the following commutative diagram
 %   \newsavebox{\tempdiagramm}
 %       \begin{lrbox}{\tempdiagramm}
        \[\begin{tikzcd}
			f_*(\omega_f) \ar[r,"\text{coev}"]
			& {[}f_*(\unitC),f_*(\omega_f)\otimes f_*(\unitC){]} \ar[r,"{[1,\text{lax}]}"] 
			& {[}f_*(\unitC),f_*(\omega_f){]} \ar[r,"{[1,\epsilon]}"] & {[}f_*(\unitC),\unitD{]} \\
			f_*(\unitC) \ar[u,"f_*(\theta)"] \ar[r,"\text{coev}"] & {[}f_*(\unitC) , f_*(\unitC) \otimes f_*(\unitC) {]} \ar[r,"{[1,\text{lax}]}"] & {[} f_*(\unitC) , f_*(\unitC) {]} \ar[u,"{[1,f_*(\theta)]}"']
        \end{tikzcd}\] 
%        \end{lrbox}
%
%        \smallskip
%       \noindent\resizebox{\linewidth}{!}{\usebox{\tempdiagramm}}
%
%        \noindent
while the adjoint of the corresponding associative form is obtained by going along the bottom.
\end{proof}

\begin{Thm}\label{thm:main-thm}
	Let $f^*\colon\cat D \to \cat C$ be a geometric functor between rigidly-compactly generated tensor-triangulated categories. Then $f^*$ is a finite \'{e}tale morphism $($Definition~\ref{def:finite-etale}\,$)$ if and only if the following three conditions hold:
	\begin{enumerate}
		\item \label{it:m-a} $f^*$ satisfies Grothendieck--Neeman duality; 
		\item \label{it:m-b} the right adjoint $f_*$ is conservative; 
		\item \label{it:m-c} the map $\unitC \to \omega_f$ adjoint to the trace map is an isomorphism.
	\end{enumerate}
\end{Thm}

\begin{proof}
	($\Rightarrow$) If $f^*$ is finite \'{e}tale, then it is extension-of-scalars with respect to the compact separable commutative algebra $f_*(\unit)$. By the separable Neeman--Thomason Localization Theorem established by Balmer \cite[Theorem~4.2]{Balmer16a}, the compact objects in $\cat C$ are precisely the thick subcategory generated by the image $f^*(\cat D^c)$ of the compact objects in $\cat D$. Thus, by the projection formula, the fact that $f_*(\unit)$ is compact ensures that $f_*(c)$ is compact for all $c \in \cat C^c$. Thus, $f^*$ satisfies Grothendieck--Neeman duality. The right adjoint $f_*$ is certainly conservative (in fact faithful). Moreover, by Corollary~\ref{cor:rigid-comm-alg}, the commutative rigid separable algebra $f_*(\unit)$ is strongly separable, so that its trace form is nondegenerate. By Lemma~\ref{lem:omega-form}, this means that the canonical map $\unitC \to \omega_f$ becomes an isomorphism after applying $f_*$. But this means the canonical map is an isomorphism since $f_*$ is conservative.

	($\Leftarrow$) If $f^*$ satisfies Grothendieck--Neeman duality, then the commutative algebra $f_*(\unitC)$ is rigid, hence has a trace map (so that part \hyperref[it:m-c]{(c)} makes sense). Moreover, by Lemma~\ref{lem:omega-form}, if the map $\unit \to \omega_f$ adjoint to the trace map is an isomorphism, then the trace form is nondegenerate; hence by Corollary~\ref{cor:rigid-comm-alg}, $f_*(\unitC)$ is a (strongly) separable algebra. By Proposition~\ref{prop:tensor-monadic}, we have a tensor-triangulated equivalence $\cat C \to f_*(\unit)\MMod_{\cat D}$ compatible with the two adjunctions. Here we use the assumption that $f_*$ is conservative and the fact that $\cat C$ is idempotent complete (since it has small coproducts). Therefore, $f^*$ is finite \'{e}tale.
\end{proof}

\begin{Rem}
	Although in part \hyperref[it:m-b]{(b)} of Theorem~\ref{thm:main-thm} we just assume $f_*$ is conservative, it follows from the other hypotheses that it is actually faithful. It also follows from \hyperref[it:m-a]{(a)} and \hyperref[it:m-c]{(c)} that $f^*$ has the full Wirthm\"{u}ller isomorphism of \cite[Theorem~1.9]{BalmerDellAmbrogioSanders16}.
\end{Rem}

\begin{Rem}
	The conservativity condition \hyperref[it:m-b]{(b)} of Theorem~\ref{thm:main-thm} can be removed if we assume instead an additional hypothesis on the category $\cat C$. The remainder of this section is devoted to explaining this modification; see Corollary~\ref{cor:monogenic-char} below.
\end{Rem}

\begin{Def}\label{def:monogenic}
	We say that a rigidly-compactly generated tensor-triangulated category $\cat T$ is \emph{monogenic} if it is generated by its unit: $\cat T = \Loc\langle \unit\rangle$. We say that~$\cat T$ is \emph{locally monogenic} if the local category $\cat T_{\cat P} \coloneqq \cat T/\Loc\langle\cat P\rangle$ is monogenic for each $\cat P \in \Spc(\cat T^c)$.
\end{Def}

\begin{Exa}\label{exa:monogenic}
	The derived category $\Der(A)$ of any commutative ring $A$ is monogenic. For any quasi-compact and quasi-separated scheme $X$, the derived category $\Derqc(X)$ is locally monogenic. Indeed, under the identification $\Spc(\Derqc(X)^c)\cong X$, prime ideals $\cat P \in \Spc(\Derqc(X)^c)$ correspond to points $x \in X$, and $\Derqc(X)_{\cat P} \cong \Der(\mathcal O_{X,x})$ is the derived category of the local ring at $x$. (See \cite[Section 4]{Balmer20} for further discussion.)
\end{Exa}

\begin{Rem}\label{rem:fully-faithful}
	A geometric functor $f^*\colon\cat D \to \cat C$ is fully faithful if and only if the unit map $\unit_{\cat D} \to f_*(\unit_{\cat C})$ is an isomorphism. Indeed, a left adjoint is fully faithful if and only if the unit of the adjunction is a natural isomorphism, and one readily checks that the unit $\eta_d \colon d \to f_*f^*(d)$ coincides with the composite
		\[
			d \simeq d \otimes \unit_{\cat D} \to d \otimes f_*(\unit_{\cat C}) \simeq f_*(f^*(d) \otimes \unit_{\cat C}) \simeq f_*f^*(d).
		\]
	Moreover, note that in this case an object $c \in\cat C$ is in the essential image of $f^*$ if and only if the counit $\epsilon_c\colon f^*f_*(c) \to c$ is an isomorphism.
\end{Rem}

\begin{Lem}\label{lem:trace-to-counit}
	If $f^* \colon \cat D \to \cat C$ is a fully faithful geometric functor, then the map $\theta\colon \unit_{\cat C} \to \omega_f$ adjoint to the trace map is an isomorphism if and only if the counit $\epsilon_{\omega_f} \colon f^*f_*(\omega_f) \to \omega_f$ is an isomorphism.
\end{Lem}

\begin{proof}
	This follows from the commutative diagram
		\[\begin{tikzcd}
			f^*f_*(\unit_{\cat C}) \ar[d,"\epsilon","\simeq"'] \ar[r,"f^*f_*(\theta)","\simeq"'] & f^*f_*(\omega_f) \ar[d,"\epsilon"] \\
			\unit_{\cat C} \ar[r,"\theta"] & \omega_f\rlap{.}
		\end{tikzcd}\]
	The left-hand counit is an isomorphism since $\unit_{\cat C} \simeq f^*(\unit_{\cat D})$ is in the essential image of $f^*$ and $f_*(\theta)$ is an isomorphism (by Lemma~\ref{lem:omega-form} and Corollary~\ref{cor:rigid-comm-alg}) since $f_*(\unit_{\cat C}) \simeq \unit_{\cat D}$ is a rigid separable commutative algebra.
\end{proof}

\begin{Prop}\label{prop:monogenic-equiv}
	Let $f^*\colon\cat D \to \cat C$ be a geometric functor between rigidly-compactly generated tensor-triangulated categories. Suppose that the following conditions hold:
		\begin{enumerate}
			\item \label{it:p-a} the unit map $\unit_{\cat D} \to f_*(\unit_{\cat C})$ is an isomorphism;
			\item \label{it:p-b} $f^*$ satisfies Grothendieck--Neeman duality; 
			\item \label{it:p-c} the map $\unit_{\cat C} \to \omega_f$ adjoint to the trace map is an isomorphism.
		\end{enumerate}
	If $\cat C$ is locally monogenic $($Definition~\ref{def:monogenic}\,$)$, then $f^*$ is an equivalence.
\end{Prop}

\begin{proof}
	For any $t \in \cat C$, consider an exact triangle
		\begin{equation}\label{eq:counit-triangle}
			f^*f_*(t) \xrightarrow{\epsilon_t} t \to \cone(\epsilon_t) \to \Sigma f^*f_*(t).
		\end{equation}
	Hypothesis \hyperref[it:p-a]{(a)} asserts that $f^*$ is fully faithful (Remark~\ref{rem:fully-faithful}). Hence ${f_*(\cone(\epsilon_t)) = 0}$. Hypotheses \hyperref[it:p-b]{(b)} and \hyperref[it:p-c]{(c)} imply (by \cite[Theorem~3.3]{BalmerDellAmbrogioSanders16}) that there is a natural isomorphism $f^* \simeq f^!$ where $f^!$ denotes the right adjoint to $f_*$. Hence the last map in \eqref{eq:counit-triangle} vanishes and we have a splitting
		\begin{equation}\label{eq:splitting}
			 t \simeq f^*f_*(t) \oplus \cone(\epsilon_t)
		 \end{equation}
	for any $t \in \cat C$. Now let $\cat P \in \Spc(\cat C^c)$ and let $L_{\cat P}\colon\cat C \to \cat C_{\cat P}$ denote the localization to the local category at $\cat P$ (see, \textit{e.g.}, \cite[Remark~1.22--Definition~1.25]{BarthelHeardSanders21pp}). This localization has an associated idempotent triangle
		\[
			e_{\cat P} \to \unit_{\cat C} \to f_{\cat P} \to \Sigma e_{\cat P}
		\]
	in $\cat C$. The kernel of $f_*$ does not contain any nonzero ring object (since if $f_*(R)=0$ then the unit $\unit_{\cat C} \to R$ would vanish by an argument similar to the proof of Lemma~\ref{lem:trace-to-counit}). Since the right idempotent $f_{\cat P}$ is a ring object, we conclude that $f^*f_*(f_{\cat P}) \neq 0$. Now the local category $\cat C_{\cat P}$ is local in the sense of \cite[Terminology~1.11]{BarthelHeardSanders21pp}. Hence by \cite[Theorem~4.5]{Balmer10b}, the endomorphism ring $\End_{\cat C_{\cat P}}(\unit_{\cat C_{\cat P}})$ is local. The identifications $\End_{\cat C_{\cat P}}(\unit_{\cat C_{\cat P}}) \simeq \Hom_{\cat C}(\unit_{\cat C},f_{\cat P}) \simeq \End_{\cat C}(f_{\cat P})$ are isomorphisms of rings, and we conclude that $\End_{\cat C}(f_{\cat P})$ has no nontrivial idempotents. Hence the fact that $f^*f_*(f_{\cat P}) \neq 0$ implies that the splitting \eqref{eq:splitting} is trivial: $f_{\cat P} \simeq f^*f_*(f_{\cat P})$.

	We have thus established that $f_{\cat P}$ (and all its shifts $\Sigma^* f_{\cat P}$) is contained in the essential image of $f^*$. By invoking $f^* \simeq f^!$ again, we conclude that if $f_*(t)=0$ then $\Hom_{\cat C}^*(t,f_{\cat P}) = 0$ for all $\cat P \in \Spc(\cat C^c)$. In particular, if $x \in \cat C^c$ is a compact(=rigid) object with $f_*(x)=0$ then
		\begin{equation}\label{eq:vanishing}
			 \Hom_{\cat C_{\cat P}}^*(L_{\cat P}(x),\unit_{\cat C_{\cat P}}) \simeq \Hom_{\cat C}^*(x,f_{\cat P}) = 0.
		\end{equation}
	Since $\cat C_{\cat P}$ is monogenic (by hypothesis) and $L_{\cat P}(x)$ is compact, \eqref{eq:vanishing} implies that $L_{\cat P}(x)=0$, so that $x \in \cat P$. This is true for all $\cat P$ (\textit{i.e.}~$\supp(x)=\emptyset$), so $x=0$. Thus, the kernel of $f_*$ does not contain any nonzero compact objects. Hence \eqref{eq:counit-triangle} shows that $f^*f_*(x)\simeq x$ for every compact $x\in \cat C^c$. It follows that $f_*$ is conservative. Indeed, if $f_*(t) = 0$, then $\Hom_{\cat C}(x,t) \simeq\Hom_{\cat C}(f^*f_*(x),t) \simeq\Hom_{\cat D}(f_*(x),f_*(t)) =0$ for every compact $x \in \cat C^c$ and hence $t=0$. Thus, since the kernel of $f_*$ is trivial, the exact triangle \eqref{eq:counit-triangle} shows that $f^*f_*(t)\simeq t$ for every $t \in \cat C$. Hence $f^*$ is essentially surjective, and the proof is complete.
\end{proof}

\begin{Exa}
	Hypothesis \hyperref[it:p-c]{(c)} in Proposition~\ref{prop:monogenic-equiv} does not follow from the other hypotheses. Indeed, consider the projective line $\Pone$ over the field $k$. The structure map $f\colon\Pone \to \Spec(k)$ induces a fully faithful geometric functor ${f^*\colon\Der(k)\to\Derqc(\Pone)}$ which satisfies Grothendieck--Neeman duality (see \cite[Example~6.14]{BalmerDellAmbrogioSanders16} or \cite{LipmanNeeman07}). Moreover, the target category is locally monogenic (see Example~\ref{exa:monogenic}). However, $\unit \not\simeq \omega_f$. Indeed, $\omega_f \simeq \Sigma \mathcal O_{\Pone}(-2)$ but
		\[
			\Hom(\Sigma \mathcal O_{\Pone}(-2),\mathcal O_{\Pone}) =\Hom(\mathcal O_{\Pone},\Sigma^{-1}\mathcal O_{\Pone}(2)) = H^{-1}(\Pone,\mathcal O_{\Pone}(2)) = 0
		\]
	while $\Hom(\mathcal O_{\Pone},\mathcal O_{\Pone}) = k$. On the other hand, since $f_*(\unit) \simeq H^*(\Pone,\mathcal O_{\Pone})=k= \unit$ is a rigid separable commutative algebra, we know from Lemma~\ref{lem:omega-form} that $f_*(\unit) \simeq f_*(\omega_f)$. Indeed, $f_*(\omega_f) \simeq H^{*+1}(\Pone,\mathcal O_{\Pone}(-2)) = k$ (see, \textit{e.g.}, \cite[Theorem 5.1]{Hartshorne77}). Note that this example is the $\cat T=\Derqc(\Pone)$ case of Example~\ref{exa:loc}.
\end{Exa}

\begin{Cor}\label{cor:monogenic-char}
	Let $f^*\colon\cat D \to \cat C$ be a geometric functor between rigidly-compactly generated tensor-triangulated categories. Suppose $\cat C$ is locally monogenic $($Definition\,\ref{def:monogenic}$)$. Then $f^*$ is a finite \'{e}tale morphism if and only if the following two conditions hold:
		\begin{enumerate}
			\item \label{it:c-a} $f^*$ satisfies Grothendieck--Neeman duality; 
			\item \label{it:c-b} the map $\unit_{\cat C} \to \omega_f$ adjoint to the trace map is an isomorphism.
		\end{enumerate}
\end{Cor}

\begin{proof}
	The $(\Rightarrow)$ direction is provided by Theorem~\ref{thm:main-thm}. We need to prove that a geometric functor $f^*\colon\cat D \to \cat C$ satisfying \hyperref[it:c-a]{(a)} and \hyperref[it:c-b]{(b)} is finite \'{e}tale provided that $\cat C$ is locally monogenic. As explained in the proof of Theorem~\ref{thm:main-thm}, hypotheses \hyperref[it:c-a]{(a)} and \hyperref[it:c-b]{(b)} imply that $f_*(\unit_{\cat C})$ is a rigid separable commutative algebra. Thus \mbox{$\cat D' \coloneqq (f_*(\unit_{\cat C})\FFree_{\cat D})^\natural \cong f_*(\unit_{\cat C})\MMod_{\cat D}$} is a rigidly-compactly generated tensor-triangulated category (Example~\ref{exa:A-modules-triangulated} and Remark~\ref{rem:AModisrcg}), and we can factor $f^*$ as a composite 
		\[
			\cat D \xrightarrow{g^*} \cat D' \xrightarrow{h^*} \cat C, 
		\]
	where $g^*$ is finite \'{e}tale and $h^*$ is a geometric functor with the property that $\unit_{\cat D'} \to h_*(\unit_{\cat C})$ is an isomorphism. (One can verify that the functor $h^*$ preserves coproducts in a routine manner using the basic properties of the Kleisli adjunction and idempotent completion. That these are symmetric monoidal functors is explained in the proof of Proposition~\ref{prop:tensor-monadic}.) Now consider a compact object $x \in \cat C^c$. Since $g^*$ is finite \'{e}tale, the counit of the $g^* \dashv g_*$ adjunction has a section (\textit{cf.}~\cite[Proposition~3.11]{Balmer11} and \cite[Theorem~1.2]{Rafael90}). Thus $h_*(x)$ is a direct summand of $g^*g_*h_*(x) = g^*f_*(x)$ and hence is compact since $f_*$ and $g^*$ preserve compactness. Thus,~$h^*$ satisfies Grothendieck--Neeman duality. Moreover, 
		\[
			\omega_h = h^!(\unit_{\cat D'}) \simeq h^!(\omega_g) = h^!(g^!(\unit_{\cat D})) \simeq f^!(\unit_{\cat D}) = \omega_f \simeq \unit_{\cat C} \simeq h^*(\unit_{\cat D'}),
		\]
	so Lemma~\ref{lem:trace-to-counit} implies that the map $\unit_{\cat C} \to \omega_h$ adjoint to the trace form on $h_*(\unit_{\cat C})$ is an isomorphism. Thus~$h^*$ satisfies all the hypotheses of Proposition~\ref{prop:monogenic-equiv}. This establishes that $h^*$ is an equivalence, and the proof is complete.
\end{proof}

\begin{Rem}
	It follows from Example~\ref{exa:local-of-etale} below that if $f^*\colon\cat D \to \cat C$ is finite \'{e}tale and $\cat D$ is locally monogenic, then $\cat C$ is locally monogenic; see Corollary~\ref{cor:mon-to-mon}. Thus, when studying the finite \'{e}tale extensions of a locally monogenic category $\cat D$, there is no loss of generality in using the criterion provided by Corollary~\ref{cor:monogenic-char}.
\end{Rem}

\begin{Rem}
	We now provide an example which shows that Proposition~\ref{prop:monogenic-equiv} and Corollary~\ref{cor:monogenic-char} do not hold (in general) without the locally monogenic hypothesis.
\end{Rem}
        
\begin{Exa}\label{exa:product-cat}
  Let $\cat T$ be a nonzero rigidly-compactly generated tensor-triangulated category. The product category $\cat T \times \cat T$ can be triangulated by defining the suspension and exact triangles coordinate-wise, and we have fully faithful triangulated functors
\[ \cat T \hookrightarrow \cat T \times \cat T \hookleftarrow \cat T \]
	given by $a \mapsto (a,0)$ and $b \mapsto (0,b)$, respectively. We can turn $\cat T \times \cat T$ into a tensor-triangulated category by defining the tensor product in a $\bbZ/2$-graded fashion, 
		\[
			(a_0,a_1) \otimes (b_0,b_1) := ((a_0 \otimes b_0) \oplus (a_1 \otimes b_1), (a_0 \otimes b_1) \oplus (a_1 \otimes b_0)),
		\]
	which can be interpreted as a Day convolution on the functor category $\cat T \times \cat T =\Fun(\bbZ/2,\cat T)$. If $\cat G$ is a set of rigid-compact generators of $\cat T$, then $(\cat G\times 0) \cup (0 \times \cat G)$ is a set of rigid-compact generators for $\cat T \times \cat T$. Moreover, the unit $(\unit,0)$ of $\cat T \times \cat T$ is compact. In summary, $\cat T \times \cat T$ is a rigidly-compactly generated tensor-triangulated category, and the inclusion $a \mapsto (a,0)$ is a fully faithful coproduct-preserving tensor-triangulated functor $f^*\colon \cat T \hookrightarrow \cat T\times \cat T$. The projection $\cat T \times \cat T \to \cat T$ onto the \supth{0} coordinate is both left and right adjoint to $f^*$. It follows that $f^*$ is a geometric functor which is not an equivalence and yet satisfies the hypotheses \hyperref[it:p-a]{(a)}, \hyperref[it:p-b]{(b)}, and~\hyperref[it:p-c]{(c)} of Proposition~\ref{prop:monogenic-equiv}. Also note that the right adjoint $f_*$ is not conservative (so $f^*$ is not finite \'{e}tale) and yet $f^*$ satisfies hypotheses \hyperref[it:c-a]{(a)} and \hyperref[it:c-b]{(b)} of Corollary~\ref{cor:monogenic-char}.
\end{Exa}

\begin{Rem}
The product category $\cat T\times \cat T$ from Example~\ref{exa:product-cat} is never locally monogenic (for $\cat T$ nonzero). This follows from our abstract theorems as noted above, but of course can be seen directly. Choose any prime $\cat P \in \Spc(\cat T^c)$. Then $\cat P \times \cat P \in \Spc((\cat T \times \cat T)^c)$. If the local category $(\cat T\times \cat T)_{\cat P \times \cat P}$ at $\cat P \times \cat P$ were monogenic, then the composite
	\[\begin{tikzcd}
		\cat T \ar[d] \ar[r,hook] & \cat T\times \cat T \ar[r] & (\cat T\times \cat T)_{\cat P \times \cat P}\ar[d,"\cong"]\\%\ar[d,phantom,sloped,"\cong"]\\
		\cat T_{\cat P} \ar[rr,hook] && (\cat T_\cat P)\times(\cat T_\cat P)
	\end{tikzcd}\]
would have a conservative right adjoint, which is a contradiction since $\cat T_{\cat P} \neq 0$.
\end{Rem}

\begin{Rem}
	Example~\ref{exa:product-cat} also demonstrates that hypothesis \hyperref[it:m-b]{(b)} of Theorem~\ref{thm:main-thm} does not follow (in general) from the other two hypotheses \hyperref[it:m-a]{(a)} and \hyperref[it:m-c]{(c)}.
\end{Rem}

\section{Examples}\label{sec:examples}

We will now discuss some examples of finite \'{e}tale morphisms with an eye to future applications.

\begin{Exa}\label{exa:finite-index}
	Let $G$ be a compact Lie group. It was proved in \cite[Theorem~1.1]{BalmerDellAmbrogioSanders15} that for any finite index subgroup ${H \le G}$, the restriction functor $\SH(G) \to \SH(H)$ between equivariant stable homotopy categories is finite \'{e}tale. We can use Theorem~\ref{thm:main-thm} to improve this to an ``if and only if'' statement: 
\end{Exa}

\begin{Thm}
	Let $G$ be a compact Lie group and let $H \le G$ be a closed subgroup. The restriction functor $\res^G_H\colon\SH(G) \to \SH(H)$ is finite \'{e}tale if and only if $H$ has finite index in $G$.
\end{Thm}

\begin{proof}
	As already mentioned, the ``if'' part is \cite[Theorem~1.1]{BalmerDellAmbrogioSanders15}. For the ``only if'' part, recall that the relative dualizing object for $\res^G_H$ is the representation sphere $S^{L(H;G)}$ for the tangent $H$-representation at the coset $eH \in G/H$ (see \cite{May03} and \cite[Remark~2.16]{Sanders19}). By Theorem~\ref{thm:main-thm}, if $\res^G_H$ is finite \'{e}tale, the canonical morphism $\unit_{\SH(H)} \to S^{L(H;G)}$ is an isomorphism. Restricting to the trivial subgroup, we obtain an isomorphism $S^0 \to S^{\dim(G/H)}$ in the nonequivariant stable homotopy category~$\SH$. The dimension of (the suspension spectrum of) a sphere is recovered by rational cohomology. Hence $\dim(G/H)=0$. The compact \mbox{0-dimensional} manifold~$G/H$ is just a finite collection of points. That is, $H$ has finite index in~$G$.
\end{proof}

\begin{Exa}\label{exa:degree-n}
	Let $p_n\colon S^1 \to S^1$ denote the degree $n$ map $z \mapsto z^n$ on the unit circle. The induced functor $p_n^*\colon \SH(S^1) \to \SH(S^1)$ is \emph{not} finite \'{e}tale (for $n \ge 2$). Indeed, this amounts to the question of whether the quotient $S^1 \to S^1/C_n$ by the subgroup of \supth{$n$} roots of unity induces a finite \'{e}tale morphism $\SH(S^1/C_n)\to\SH(S^1)$. But \cite[Proposition~3.2]{Sanders19} establishes that inflation $\infl_{G/N}^G \colon  \SH(G/N) \to \SH(G)$ never satisfies Grothendieck--Neeman duality except when $N=1$ is the trivial subgroup.
\end{Exa}

\begin{Rem}
	Another way of appreciating why Example~\ref{exa:degree-n} is not finite \'{e}tale is to look at its behaviour on the Balmer spectrum, which we know due to \cite{BarthelGreenleesHausmann20,BalmerSanders17}. The points of $\Spc(\SH(S^1)^c)$ are of the form $\cat P(H,\cat C)$ for $H$ a closed subgroup of $S^1$ and $\cat C \in \Spc(\SH^c)$. The closed subgroups of~$S^1$ are, in addition to~$S^1$ itself, the finite cyclic groups $C_m$ ($m \ge 1$) realized as the roots of unity in~$S^1$. Consider the map on the Balmer spectrum
	\[
		\varphi \coloneqq \Spc(p_n^*) \colon \Spc(\SH(S^1)^c) \to \Spc(\SH(S^1)^c)
	\]
	induced by the degree $n$ map $p_n\colon S^1 \to S^1$. One can show that $\varphi(\cat P(C_m,\cat C)) = \cat P(C_{\lcm(m,n)/n},\cat C)$. For example, taking $n=2$ and fixing the nonequivariant prime~$\cat C$, it maps
	\[
		m \mapsto \begin{cases}
		m/2 & \text{if } 2\mid m,\\
		m & \text{if } 2 \nmid m.
		\end{cases}
	\]
	In particular, we find that the fibers have cardinality
	\[
		|\varphi^{-1}(\{\cat P(C_N,\cat C)\})| = \begin{cases} 1 &\text{if } 2 \mid N,\\
		2 & \text{if } 2\nmid N.
		\end{cases}
	\]
	For example, the fiber over $\cat P(C_1,\cat C)$ consists of two points:  $\big\{ \cat P(C_1,\cat C), \cat P(C_2,\cat C) \big\}$. Moreover, if the nonequivariant prime~$\cat C\coloneqq \cat C_{2,\infty}$ is the $2$-local prime at chromatic height $\infty$, then $\cat P(C_1,\cat C_{2,\infty}) \subseteq \cat P(C_2,\cat C_{2,\infty})$ is a nontrivial inclusion in the fiber over $\cat P(C_1,\cat C_{2,\infty})$. This implies that the basic theorems of Balmer \cite[Theorem~1.5]{Balmer16b} on the behaviour of finite \'{e}tale morphisms do not hold for the morphisms $p_n^*\colon \SH(S^1)\to\SH(S^1)$.
\end{Rem}

\begin{Lem}\label{lem:beck}
	Consider a diagram of coproduct-preserving $(N$-$)$tensor-triangulated functors between rigidly-compactly generated $(N$-$)$tensor-triangulated categories
	\[
		\begin{tikzcd}
			\cat C \ar[r,"g^*"] \ar[d,"h^*"'] & \cat D \ar[d,"k^*"] \\
			\cat C' \ar[r,"f^*"'] & \cat D'
		\end{tikzcd}
	\]
	which commutes up to natural isomorphism of symmetric monoidal functors. Denote the right adjoints by $f^* \dashv f_*$ and $g^* \dashv g_*$, and suppose that the Beck--Chevalley comparison map
	\[
		h^*g_* \to f_*k^*
	\]
	is a natural isomorphism. If $g^*$ is finite \'{e}tale and $f_*$ is conservative, then $f^*$ is finite \'{e}tale.
\end{Lem}

\begin{proof}
	The Beck--Chevalley comparison map $h^*g_* \to f_* k^*$ is a monoidal natural transformation between lax symmetric monoidal functors; hence the natural isomorphism $h^*g_*\isor f_*k^*$ provides an isomorphism of commutative algebras $h^*g_*(\unitD) \simeq f_*(\unit_{\cat D'})$. By assumption, $g_*(\unitD)$ is a compact commutative separable algebra in $\cat C$; hence $f_*(\unit_{\cat D'})$ is a compact commutative separable algebra in~$\cat C'$. The $f^* \dashv f_*$ adjunction satisfies the projection formula (see \cite[Proposition~2.15]{BalmerDellAmbrogioSanders16}), and $f_*$ is conservative by hypothesis. Hence, Proposition~\ref{prop:tensor-monadic} provides the result.
\end{proof}

\begin{Exa}\label{exa:A-modules}
	If $\cat C$ is a presentably symmetric monoidal stable $\infty$-category and $A \in \CAlg(\cat C)$ is a commutative algebra in $\cat C$, then we can consider the presentably symmetric monoidal stable $\infty$-category $A\MMod_{\cat C}$ of $A$-modules. If $\cat C$ is rigidly-compactly generated, then so is $A\MMod_{\cat C}$ (see \cite[Remark~3.11]{PatchkoriaSandersWimmer22}, for example). At the level of homotopy categories, the extension-of-scalars $\Ho(\cat C) \to \Ho(A\MMod_{\cat C})$ is then a geometric functor of rigidly-compactly generated $\infty$-tensor-triangulated categories whose right adjoint is conservative.
\end{Exa}

\begin{Exa}\label{exa:infty-cat-algebra}
	Let $\cat C$ be a presentably symmetric monoidal stable $\infty$-category, and let $A,B \in \CAlg(\cat C)$ be commutative algebras in $\cat C$. We then have
	\[
		\begin{tikzcd}
			\cat C \ar[r]  \ar[d] & B\MMod_{\cat C} \ar[d] \\
			A\MMod_{\cat C} \ar[r] & (A \otimes B)\MMod_{\cat C}\rlap{,}
		\end{tikzcd}
	\]
	where all four functors are extension-of-scalars. This is an example where the Beck--Chevalley property holds (at the level of the underlying stable $\infty$-categories). In particular, the induced diagram of $\infty$-tensor-triangulated categories
	\[
		\begin{tikzcd}
			\Ho(\cat C) \ar[r]  \ar[d] & \Ho(B\MMod_{\cat C}) \ar[d] \\
			\Ho(A\MMod_{\cat C}) \ar[r] & \Ho((A \otimes B)\MMod_{\cat C})
		\end{tikzcd}
	\]
	satisfies the first hypothesis of Lemma~\ref{lem:beck}. Moreover, the right adjoints are all conservative (Example~\ref{exa:A-modules}). Thus, if the top horizontal functor is finite \'{e}tale (\textit{i.e.}~if~$B$ is a compact separable commutative algebra in $\Ho(\cat C)$), then the bottom horizontal functor is also finite \'{e}tale.
\end{Exa}

\begin{Exa}
	Let $G$ be a compact Lie group and let $\Sp_G$ denote the symmetric monoidal stable \mbox{$\infty$-category} of $G$-spectra (see \cite[Appendix C]{GepnerMeier}). Let $\triv_G \colon \Sp \to \Sp_G$ denote the unique colimit-preserving symmetric monoidal functor from the $\infty$-category of spectra. Since $\res^G_H \circ \triv_G \simeq \triv_H$ for any $H \le G$, we have a commutative diagram
	\[
		\begin{tikzcd}
			\Ho(\Sp_G) \ar[d] \ar[r] & \Ho(\Sp_H) \ar[d] \\	
			\Ho(\triv_G(\bbE)\MMod_{\Sp_G}) \ar[r] & \Ho(\triv_H(\bbE)\MMod_{\Sp_H})
		\end{tikzcd}
	\]
	for any $\bbE \in \CAlg(\Sp)$. If $H \le G$ has finite index, then the top horizontal functor is finite \'{e}tale (Example~\ref{exa:finite-index}), and hence the bottom horizontal functor is finite \'{e}tale. Taking $\bbE = \HZ$, we obtain that the restriction functor
	\[
		\DHZG \to \DHZH
	\]
	between categories of derived Mackey functors studied in \cite{PatchkoriaSandersWimmer22} is finite \'{e}tale. This will be utilized in the forthcoming \cite{BarthelHeardSanders21pp} which will classify the localizing tensor-ideals of these categories.
\end{Exa}

\begin{Exa}\label{exa:A-B-separable}
	A version of Example~\ref{exa:infty-cat-algebra} holds purely at the level of triangulated categories if one assumes the two algebras are separable. More precisely, let $A$ and $B$ be two commutative separable algebras in a rigidly-compactly generated $N$-tensor-triangulated category $\cat T$. Iterated extension-of-scalars behaves as one expects (see \cite[Proposition~1.14]{Pauwels17}) and we have a diagram of rigidly-compactly generated $N$-tensor-triangulated categories
	\[
		\begin{tikzcd}
			\cat T \ar[r,"F_A"] \ar[d,"F_B"] & A\MMod_{\cat T} \ar[d]\\
			B\MMod_{\cat T} \ar[r] & (A\otimes B)\MMod_{\cat T}
		\end{tikzcd}
	\]
	which commutes up to isomorphism. Lemma~\ref{lem:beck} implies that if the top functor is finite \'{e}tale, then so is the bottom functor.
\end{Exa}

\begin{Rem}\label{rem:restriction}
	Let $F\colon \cat D \to \cat C$ be a geometric functor of rigidly-compactly generated tensor-triangulated categories and let $\varphi\colon \Spc(\cat C^c) \to \Spc(\cat D^c)$ be the induced map on spectra. For any Thomason subset $Y\subseteq \Spc(\cat D^c)$ with $V\coloneqq \Spc(\cat D^c)\setminus Y$, we have an induced functor $F|_{V}\colon \cat D(V) \to \cat C(\varphi^{-1}(V))$ on finite localizations such that
	\begin{equation}\label{eq:localization-in-target}
		\begin{tikzcd}
			\cat D \ar[d] \ar[r,"F"] & \cat C \ar[d]\\
			\cat D(V) \ar[r,"F|_{V}"] & \cat C(\varphi^{-1}(V))
		\end{tikzcd}
	\end{equation}
	commutes up to isomorphism. Moreover, on spectra,  
	\[
		\varphi^{-1}(V)\cong \Spc(\cat C(\varphi^{-1}(V))) \xra{\Spc(F|_{V})} \Spc(\cat D(V))\cong V
	\]
	is just the restriction $\varphi|_{V} \colon \varphi^{-1}(V) \to V$.
\end{Rem}

\begin{Exa}[Restriction in the target]\label{exa:local-of-etale}
	If $F\colon \cat D \to \cat C$ is finite \'{e}tale, then the induced ``restriction'' functor
	\[
		F|_{V} \colon \cat D(V) \to \cat C(\varphi^{-1}(V))
	\]
	of Remark~\ref{rem:restriction} is also finite \'{e}tale. Here $V \subseteq \Spc(\cat D^c)$ is the complement of a Thomason subset. For example,~$V$ could be a quasi-compact open subset. Indeed, this is just a special case of Example~\ref{exa:A-B-separable} with $B=f_{V^c}$ the idempotent algebra for the finite localization $\cat D \to \cat D(V)$.
\end{Exa}

\begin{Cor}\label{cor:mon-to-mon}
	Let $F\colon \cat D \to \cat C$ be finite \'{e}tale. If $\,\cat D$ is locally monogenic, then $\cat C$ is locally monogenic.
\end{Cor}

\begin{proof}
	We use the notation of Remark~\ref{rem:restriction}. Let $\cat P \in \Spc(\cat C^c)$ and consider its image $\varphi(\cat P) \in \Spc(\cat D^c)$. Let $V\coloneqq \gen(\varphi(\cat P))$ be the subset of $\Spc(\cat D^c)$ consisting of all generalizations of the point $\varphi(\cat P)$. It is the complement of a Thomason subset, and restriction to $V$ is localization at the point $\varphi(\cat P)$; see \cite[Remark~1.21 and Definition~1.25]{BarthelHeardSanders21pp}. Note that $\gen(\cat P) \subseteq \varphi^{-1}(V)$. Then consider the composition
	\[
		\cat D_{\varphi(\cat P)} = \cat D(V) \to \cat C(\varphi^{-1}(V)) \to \cat C(\gen(\cat P))=\cat C_{\cat P}, 
	\]
	where the first functor is finite \'{e}tale (Example~\ref{exa:local-of-etale}) and the second functor is a localization. Since the right adjoints are conservative, any set of compact generators of $\cat D_{\varphi(\cat P)}$ is mapped to a set of compact generators of $\cat C_{\cat P}$. Thus, $\cat D_{\varphi(\cat P)}$ monogenic implies $\cat C_{\cat P}$ monogenic.
\end{proof}

\begin{Rem}
	Additional equivariant examples are featured in the work of Balmer and Dell'Ambrogio on Mackey 2-motives \cite{BalmerDellAmbrogio20,DellAmbrogio22}. On the other hand, the following basic example relates the tensor-triangular notion of finite \'{e}tale with the ordinary scheme-theoretic notion: 
\end{Rem}

\begin{Thm}[Balmer]\label{thm:Balmer}
	If $f\colon X\to Y$ is a finite \'{e}tale morphism of quasi-compact and quasi-separated schemes, then the derived functor $\bbL f^*\colon \Derqc(Y) \to \Derqc(X)$ is a finite \'{e}tale morphism in the sense of Definition~\ref{def:finite-etale}.
\end{Thm}

\begin{proof}
	This is provided by \cite[Theorem~3.5]{Balmer16a}; see also \cite[Example 0.3]{Neeman18}.
\end{proof}

\begin{Rem}
	The proof of the above theorem works verbatim for other tensor-triangulated categories $\cat T(X)$ fibered over a category of schemes, provided the pseudofunctor $X \mapsto \cat T(X)$ satisfies flat base change. Many motivic examples of such pseudofunctors are discussed in \cite{CisinskiDeglise19}. We just mention:
\end{Rem}

\begin{Exa}
	Let $\mathbb{L}/\mathbb{K}$ be a finite separable extension of fields whose characteristic (if positive) is invertible in the ring $R$. The induced functor $\SH(\mathbb{K};R) \to \SH(\mathbb{L};R)$ between motivic stable homotopy categories (with coefficients in $R$) is a finite \'{e}tale morphism in the sense of Definition~\ref{def:finite-etale}. The same is true of the induced functor $\DM(\mathbb{K};R) \to \DM(\mathbb{L};R)$ between derived categories of motives. See \cite{CisinskiDeglise19,Ayoub07a,Ayoub07b,Totaro18} for more information about these categories. The assumption on the characteristic ensures that these categories are rigidly-compactly generated. To see this, first recall that $\SH(\mathbb{K};R)$ is compactly generated by the twists of smooth $\mathbb{K}$-schemes of finite type (see \cite[Corollary~1.3]{Riou05} and \cite[Theorem~4.5.67]{Ayoub07b}) and, similarly, $\DM(\mathbb{K};R)$ is compactly generated by the twists of smooth separated $\mathbb{K}$-schemes of finite type (see \cite[Section~11.1]{CisinskiDeglise19} and \cite[Lemma~5.4]{Totaro16}). It is nontrivial that these generators are dualizable. For the motivic stable homotopy category, see \cite[Theorem~3.2.1]{ElmantoKhan20}, which builds on \cite[Appendix~B]{LevineYangZhaoRiou19}; for the derived category of motives, see \cite[Lemma~5.5]{Totaro16}, which extends \cite[Theorem~5.5.14]{Kelly_dissertation} and \cite[Theorem.~4.3.7]{Voevodsky00}. Since the unit object is compact and there is a generating set of dualizable objects, it follows that the compact and dualizable objects coincide. This follows from \cite[Theorem~A.2.5.(a)]{HoveyPalmieriStrickland97} and \cite[Lemma~2.2]{Neeman92b} (which is also proved in \cite[Proposition~2.1.24]{Ayoub07a}).
\end{Exa}

\begin{Rem}
	The author thinks it is interesting to have an ``intrinsic'' characterization of finite \'{e}tale morphisms in tensor triangular geometry as expressed in Theorem~\ref{thm:main-thm}. Nevertheless, actually classifying the finite \'{e}tale extensions of a given category $\cat T$ amounts to classifying the rigid (strongly) separable commutative algebras in $\cat T$. For the equivariant stable homotopy category $\cat T=\SH(G)$, this classification will be studied in forthcoming work with Balmer. The analogous problem for the stable module category $\cat T=\StMod(kG)$ has been studied in \cite{BalmerCarlson18} and is surprisingly subtle. It is currently only understood when $G$ is cyclic.
\end{Rem}

\begin{Rem}
	For the derived category $\cat T=\Derqc(X)$ of a noetherian scheme, Neeman~\cite{Neeman18} has obtained a very satisfactory classification of the (not necessarily compact) commutative separable algebras. His work shows that the tensor-triangular analogue of \'{e}tale morphism (a.k.a.~extension by a commutative separable algebra) lies somewhere between the classical \'{e}tale morphisms of schemes and the pro-\'{e}tale morphisms of Bhatt--Scholze \cite{BhattScholze15}. His results also show that there are no exotic \'{e}tale extensions of derived categories of schemes: An \'{e}tale extension of a derived category of a scheme is another derived category of a scheme. We will state this result precisely in the case of finite \'{e}tale extensions: 
\end{Rem}

\begin{Thm}[Neeman]
	Let $X$ be a noetherian scheme. If $F\colon \Derqc(X) \to \cat S$ is a finite \'{e}tale morphism $($Definition~\ref{def:finite-etale}$)$, then there exist a finite \'{e}tale morphism of schemes $f\colon U \to X$ and a tensor-triangulated equivalence $\cat S \cong \Derqc(U)$. With this identification, $F$ is naturally isomorphic to $\bbL f^*\colon \Derqc(X) \to \Derqc(U)$.
\end{Thm}

\begin{proof}
	Let $G$ denote the right adjoint of $F$. By definition, $F$ is extension-of-scalars with respect to the compact commutative separable algebra $G(\unit) \in \Derqc(X)$. Neeman \cite[Theorem~7.10]{Neeman18} establishes that there is a separated finite-type \'{e}tale map of schemes $g\colon V\to X$ and a generalization-closed subset $U \subset V$ such that $G(\unit) \cong \bbR f_*(\cat O_U)$, where $f\colon U\to X$ denotes the composite $U \hookrightarrow V \xra{g} X$. It then follows from Proposition~\ref{prop:tensor-monadic} that $\cat S \cong \Derqc(U)$ with $F \cong \bbL f^*$. Now, since $G(\unit)$ is compact, the argument in \cite[Remark~0.6]{Neeman18} shows that $U \subset V$ is actually an open subset. (Take $L\coloneqq 0$, $\widetilde{K} \coloneqq f_* f^*(K)$ and the identity map $\widetilde{K} \to f_*f^* K$ in \textit{loc.~cit.}) Thus, $f\colon U\to X$ is a separated finite-type \'{e}tale map. It is also proper since $\bbL f^* \cong F$ satisfies GN-duality (by \cite{LipmanNeeman07}; see also \cite[Section 7]{Sanders19} and \cite[Section 4.3]{Lipman09}). This completes the proof since an \'{e}tale map is proper if and only if it is finite.
\end{proof}

\begin{Rem}
	For the purpose of classifying the \emph{finite} \'{e}tale extensions of a given tensor-triangulated category, the results of Section~\ref{sec:strongly-separable} are worth keeping in mind. They clarify that the compact/rigid commutative separable algebras that provide finite \'{e}tale extensions are necessarily self-dual. This puts limits on the role finite \'{e}tale morphisms can play in equivariant contexts over non-finite groups. Stated differently, Theorem~\ref{thm:main-thm} shows that the relative dualizing object $\omega_f$ for a finite \'{e}tale morphism $f^*$ must be trivial. It is natural to wonder if there is a reasonable generalization of ``finite \'{e}tale'' in tensor triangular geometry which shares some of its good properties (\textit{e.g.}, the results of \cite{Balmer16a,Balmer16b}) and yet covers examples having non-trivial dualizing objects (\textit{e.g.}, the examples which arise in \cite{Rognes08}). 
\end{Rem}

%%%%%%%%%%%%%%%%%%%%%
% References
%%%%%%%%%%%%%%%%%%%%%

\newcommand{\etalchar}[1]{$^{#1}$}

\end{document}